\magnification=1200

\loadmsam
\loadmsbm
\loadeufm
\loadeusm
\UseAMSsymbols

\font\BIGtitle=cmr10 scaled\magstep3
\font\bigtitle=cmr10 scaled\magstep1
\font\boldsectionfont=cmb10 scaled\magstep1
\font\section=cmsy10 scaled\magstep1

\def\scr#1{{\fam\eusmfam\relax#1}}

\def\scrH{{\scr H}}

\def\scrJ{{\scr J}}

\def\scrM{{\scr M}}

\def\scrV{{\scr V}}
\def\scrX{{\scr X}}

\def\gr#1{{\fam\eufmfam\relax#1}}

	\def\grg{{\gr g}}
	\def\grh{{\gr h}}

	\def\grk{{\gr k}}
\def\grL{{\gr L}}

\def\db#1{{\fam\msbfam\relax#1}}

\def\dbA{{\db A}} 
\def\dbC{{\db C}} 
 \def\dbF{{\db F}}
\def\dbG{{\db G}}

 \def\dbN{{\db N}}
 
\def\dbQ{{\db Q}} \def\dbR{{\db R}}

 \def\dbZ{{\db Z}}

\def\Ker{\text{Ker}}

\def\der{\text{der}}
\def\d{\text{d}}

\def\sc{\text{sc}}
\def\s{\text{s}}
\def\Res{\text{Res}}
\def\ab{\text{ab}}
\def\ad{\text{ad}}

\def\Gal{\text{Gal}}
\def\Hom{\text{Hom}}
\def\End{\text{End}}

\def\Spec{\text{Spec}}

\def\Lie{\text{Lie}}

\def\leaderfill{\leaders\hbox to 1em
     {\hss.\hss}\hfill}
\def\nspace{\lineskip=1pt\baselineskip=12pt\lineskiplimit=0pt}

\def\finishproclaim{\par\rm
     \ifdim\lastskip<\medskipamount\removelastskip
     \penalty55\medskip\fi}
\def\endproof{$\hfill \square$}
\def\proof{\par\noindent {\it Proof:}\enspace}
\def\references#1{\par
  \centerline{\boldsectionfont References}%
     \parindent=#1pt\nspace}
\def\Ref[#1]{\par\hang\indent\llap{\hbox to\parindent
     {[#1]\hfil\enspace}}\ignorespaces}
\def\Item#1{\par\hang\indent\llap{\hbox to\parindent
     {#1\hfill$\,\,$}}\ignorespaces}
\def\ItemItem#1{\par\indent\hangindent2\parindent
     \hbox to \parindent{#1\hfill\enspace}\ignorespaces}

\def\Le{{\mathchoice{\,{\scriptstyle\le}\,}
  {\,{\scriptstyle\le}\,}
  {\,{\scriptscriptstyle\le}\,}{\,{\scriptscriptstyle\le}\,}}}
\def\Ge{{\mathchoice{\,{\scriptstyle\ge}\,}
  {\,{\scriptstyle\ge}\,}
  {\,{\scriptscriptstyle\ge}\,}{\,{\scriptscriptstyle\ge}\,}}}

\def\arrowsim{\,\smash{\mathop{\to}\limits^{\lower1.5pt
  \hbox{$\scriptstyle\sim$}}}\,}

\def\doublemaprights#1#2#3#4{\raise3pt\hbox{$\mathop{\,\,\hbox to
     #1pt{\rightarrowfill}\kern-30pt\lower3.95pt\hbox to
     #2pt{\rightarrowfill}\,\,}\limits_{#3}^{#4}$}}

\def\rightcapdownarrow{\raise9pt\hbox{$\ssize\cap$}\kern-7.75pt
     \Big\downarrow}

\def\rcapmapdown#1{\rightcapdownarrow\kern-1.0pt\vcenter{
     \hbox{$\scriptstyle#1$}}}

\def\rmapdown#1{\Big\downarrow\kern-1.0pt\vcenter{
     \hbox{$\scriptstyle#1$}}}
\def\rightsubsetarrow#1{{\ssize\subset}\kern-4.5pt\lowe r2.85pt
     \hbox to #1pt{\rightarrowfill}}
\def\longtwoheadedrightarrow#1{\raise2.2pt\hbox to #1pt{\hrulefill}
     \!\!\!\twoheadrightarrow}

\def\Gal{\operatorname{\hbox{Gal}}}
\def\Hom{\operatorname{\hbox{Hom}}}

\def\im{\hbox{Im}}

\NoBlackBoxes
\parindent=25pt
\document
\footline={\hfil}

\null
\vskip 0.25 in
\noindent
{\bigtitle Adrian Vasiu}
\vskip 0.15in
\noindent
{\BIGtitle Surjectivity Criteria for p-adic Representations, Part II}
\footline={\hfill}
\vskip 0.2in
\noindent
{Received: 25.02.2003}
\noindent
{/ Revised versions: 25.07.2003, 26.01.2004 and 23.03.2004}
\vskip 0.2in
\noindent
{\bf ABSTRACT}. We prove surjectivity criteria for $p$-adic representations and we apply them to abelian varieties over number fields. In particular, we provide examples of Jacobians over $\dbQ$ of dimension $d\in\{1,2,3\}$ whose $2$-adic representations have as images $GSp_{2d}(\dbZ_2)$. 
\noindent
$\vfootnote{} {Adrian Vasiu, Mathematics Department, University of Arizona, 617 N. Santa Rita, P.O. Box 210089, Tucson, AZ-85721, USA. e-mail: adrian\@math.arizona.edu}$
\noindent
$\vfootnote{} {{\it Mathematics Subject Classification (2000)}: Primary 11S23, 14L17, 17B45, 20G05 and 20G40}$

\footline={\hss\tenrm \folio\hss}
\pageno=1

\bigskip
\noindent
{\boldsectionfont \S1. Introduction}

\bigskip
Let $\Spec(R)$ be a connected affine scheme. If $M$ is an $R$-module, let $GL(M)$ be the group scheme over $R$ of linear automorphisms of $M$. If $*_R$ or $*$ is an object of the category of $R$-schemes, let $*_U$ be its pull back via an affine morphism $\Spec(U)\to\Spec(R)$. A reductive group scheme $F$ over $R$ has connected fibres. Let $F^{\der}$, $Z(F)$, $F^{\ab}$, and $F^{\ad}$ be the derived group, the center, the maximal abelian quotient, and respectively the adjoint group of $F$. So $F^{\ad}=F/Z(F)$ and $F^{\ab}=F/F^{\der}$. Let $Z^0(F)$ be the maximal torus of $Z(F)$. Let $F^{\sc}$ be the simply connected semisimple group cover of $F^{\der}$. Let $c(F^{\der})$ be the degree of the central isogeny $F^{\sc}\to F^{\der}$.

Let $E$ be a number field. We fix an embedding $i_E:E\hookrightarrow\dbC$ and we identify naturally $\overline{E}=\overline{\dbQ}\subset\dbC$. Let $A$ be an abelian variety over $E$ of dimension $d\in\dbN$. To ease notations, let $H_1:=H_1(A_{\dbC},\dbZ)$
be the first Betti homology group of $A_{\dbC}$ with coefficients in $\dbZ$ and let $H_A$ be the Mumford--Tate group of $A_{\dbC}$. We recall that $H_A$ is a reductive group over $\dbQ$ and that $H_A$ is the smallest subgroup of $GL(H_1\otimes_{\dbZ}\dbQ)$ such that the Hodge cocharacter $\mu_A\colon\dbG_m\to GL(H_1\otimes_{\dbZ} \dbC)$
factors through $H_{A\dbC}$, cf. [9, 2.11, 3.4 and 3.6]. If $H_1\otimes_{\dbZ} \dbC=F^{-1,0}\oplus F^{0,-1}$ is the usual Hodge decomposition, then $\beta\in\dbG_m(\dbC)$ acts through $\mu_A$ trivially on $F^{0,-1}$ and as the multiplication with $\beta$ on $F^{-1,0}$. 

Let $p\in\dbN$ be a prime. Let $T_p(A)$ be the Tate-module of $A$. As a $\dbZ_p$-module we identify it canonically with $H_1\otimes_{\dbZ} \dbZ_p$. Let $H_{A\dbZ_p}$ be the Zariski closure of $H_{A\dbQ_p}$ in $GL(T_p(A))=GL(H_1\otimes_{\dbZ} \dbZ_p)$ and let
$$\rho_{A,p}\colon\Gal(E)\to GL(T_p(A))(\dbZ_p)=GL(H_1\otimes_{\dbZ} \dbZ_p)(\dbZ_p)$$
be the $p$-adic representation. Let $V_p(A):=T_p(A)[{1\over p}]=H_1\otimes_{\dbZ} \dbQ_p$. Let $G_p$ be the connected algebraic subgroup of $GL(V_p(A))=GL(H_1\otimes_{\dbZ} \dbQ_p)$ that is the identity component of the Zariski closure of $\im(\rho_{A,p})$ in $GL(V_p(A))$. It is known that there is a smallest finite field extension $E^{\text{conn}}$ of $E$ such that the compact group
$$K_p(A):=\rho_{A,p}(\Gal(E^{\text{conn}}))$$ 
is an open subgroup of $G_p(\dbQ_p)$ (see [4]). The field $E^{\text{conn}}$ does not depend on $p$, cf. [33, p. 15] or [20, Th. 0.1]. A theorem obtained independently by Deligne, Borovoi and Pyatetskii-Shapiro asserts that $G_p$ is a subgroup of $H_{A\dbQ_p}$ (for instance, cf. [9, 2.9 and 2.11]). So $K_p(A)$ is a subgroup of $H_{A\dbZ_p}(\dbZ_p)$. 

\bigskip\noindent
{\bf 1.1. The problem.} The problem we deal with is to describe the subgroup $K_p(A)$ of $H_{A\dbZ_p}(\dbZ_p)$. It has a long history for elliptic curves (see [21], [27], etc.) and for $p>>0$ (see [33], [19], [24], [11], [42], etc.). For instance, if $A$ is a semistable elliptic curve over $\dbQ$ (so $d=1$ and $E=\dbQ$), then Serre proved that for $p\Ge 3$ the homomorphism $\rho_{A,p}$ is surjective iff its reduction mod $p$ is irreducible (cf. [27, Prop. 21], [28, p. IV-23-24] and [31, \S3.1]). In this Part II we work in a context in which the following three conditions hold:

\medskip
{\bf (a)} the Mumford--Tate conjecture holds for $(A,p)$ (i.e. we have $G_p=H_{A\dbQ_p}$);

\smallskip
{\bf (b)} the group $H_{A\dbZ_p}$ is a reductive group scheme over $\dbZ_p$; 

\smallskip
{\bf (c)} the group $K_p(A)$ surjects onto $H_{A\dbZ_p}/Z^0(H_{A\dbZ_p})(\dbF_p)$.

\medskip
Fixing $A$, it is known that (a) implies that for $p>>0$ we have $\im(H_{A\dbZ_p}^{\sc}(\dbZ_p)\to H_{A\dbZ_p}^{\der}(\dbZ_p))\vartriangleleft K_p(A)$ (see [42, Th. 2]; see also [19]). The main goal of this Part II is to get several higher dimensional analogues of Serre's result as well as some concrete forms of Wintenberger and Larsen's results. The refined study of Part I (see [37]) allows us to treat also the harder cases when either $p$ is small (like $p=2$ or $p=3$) or we are in some exceptional situation (like when $p$ divides either $c(H_A^{\der})$ or the order of $Z(H_A^{\der})$). Condition (a) is almost always implied by (b) and (c). Condition (b) is often implied by different geometric assumptions on $A$ and $\End_{\overline{E}}(A)$. So in \S4 we gather examples stated without assuming that conditions (a) to (c) hold. In particular, in 4.1.1 and 4.2 we use Jacobians over $\dbQ$ of dimension $d\in\{1,2,3\}$ to provide concrete situations when the groups $K_2(A)$, $\im(\rho_{A,2})$ and $H_{A\dbZ_2}(\dbZ_2)$ are all equal and isomorphic to $GSp_{2d}(\dbZ_2)$. We now describe the contents of \S2 and \S3.

\bigskip\noindent
{\bf 1.2. Abstract theory.} In \S2 we work with a reductive group scheme $G$ over the Witt ring $W(k)$ of a finite field $k$ and with a closed subgroup $K$ of $G(W(k))$. The results [37, 1.3, 4.1.1, 4.5] implicitly classify all semisimple group schemes $G$ such that $G(W(k))$ has only one closed subgroup surjecting onto $G(k)$ (see 2.2.5). But for a refined study of $K_p(A)$ this classification does not suffice. So in 2.3 and 2.4 we do not assume that $Z^0(G)$ is nontrivial. If $\im(K\to G(k))$ contains $G^{\der}(k)$ and if $g.c.d.(p,c(G^{\der}))=1$, then outside a precise (finite) list of exceptions we show that $G^{\der}(W(k))\vartriangleleft K$ (see 2.3). But this does not imply that $K$ surjects onto $(G/Z^0(G))(W(k))$. Theorem 2.4.1  identifies some conditions that imply that $K$ surjects onto $(G/Z^0(G))(W(k))$. In particular, Theorem 2.4.1 applies if $G^{\der}$ is a form of $SL_{p^2n}/\mu_p$ or if $p=2$, $c(G^{\der})=2$ and $G/Z^0(G)$ is adjoint of either $D_{2n+1}$ or ${}^2D_{2n+1}$ Lie type. See 2.4.2 and 2.4.3 for variants and examples to 2.4.1.

\bigskip\noindent
{\bf 1.3. Geometric theory.} In \S 3 we apply \S 2 to study $K_p(A)$ under the conditions 1.1 (a) to (c). It is expected that under a certain ``maximal" assumption on the factorization of $\mu_A$ through $H_{A\dbC}$, we have $K_p(A)=H_{A\dbZ_p}(\dbZ_p)$ for $p>>0$ (see [30, \S 9 to \S11]). So if this ``maximal" condition holds (resp. does not hold), then it is of interest to recognize when the equality $K_p(A)=H_{A\dbZ_p}(\dbZ_p)$ (resp. the inclusion $\im(H_{A\dbZ_p}^{\sc}(\dbZ_p)\to H_{A\dbZ_p}^{\der}(\dbZ_p))\vartriangleleft K_p(A)$) holds, by working if possible only mod $p$. 

For the study of $K_p(A)$ in \S3 we either assume that $H_{A\dbZ_p}^{\ab}$ is $\dbG_m$ (and so we can ``appeal" to the assumed to be surjective $p$-adic cyclotomic character of $\Gal(\overline{E}/E^{\text{conn}})$) or take one of the following two approaches (see 3.3 and 3.4). In the first one we state results in terms of the index $[H_{A\dbZ_p}(\dbZ_p):K_p(A)]$. In the second one we get surjectivity criteria onto $H_{A\dbZ_p}/Z^0(H_{A\dbZ_p})(\dbZ_p)$. If $g.c.d.(p,c(H_{A\dbZ_p}/Z^0(H_{A\dbZ_p})))=1$, then these results and criteria are in essence consequences of [37, Th. 1.3]. So we now mention two applications of 2.4.1 and 2.4.2 with $p$ dividing $c(H_{A\dbZ_p}/Z^0(H_{A\dbZ_p}))$. 

The first one is a mod $2$ surjectivity criterion onto $H_{A\dbZ_2}/Z^0(H_{A\dbZ_2})(\dbZ_2)$ in contexts with $H_{A\dbZ_2}^{\ab}=\dbG_m$ and $H_{A\dbZ_2}^{\ad}$ absolutely simple; see 3.4.1.1 (a). The second one is a mod $p$ surjectivity criterion onto $H_{A\dbZ_p}/Z^0(H_{A\dbZ_p})(\dbZ_p)$ (see 3.4.2.1) that combines our results with standard properties of $p$-adic semistable representations; it involves an arbitrary $p$ and an $H_{A\dbZ_p}$ that is a $\dbG_m\times GL_{pn}$ group.

Each criterion of \S3 has a version where instead of assuming that $K_p(A)$ surjects onto $H_{A\dbZ_p}/Z^0(H_{A\dbZ_p})(\dbF_p)$ we only assume that $\im(H_{A\dbZ_p}^{\sc}(\dbF_p)\to H_{A\dbZ_p}/Z^0(H_{A\dbZ_p})(\dbF_p))\vartriangleleft\im(K_p(A)\to H_{A\dbZ_p}/Z^0(H_{A\dbZ_p})(\dbF_p))$. Though these variants are especially important for the situations when the hinted at ``maximal" condition on $\mu_A$ does not hold, not to make this Part II too long they will not be stated here (they require no extra arguments).

                                                                               \bigskip
\noindent                                             
{\boldsectionfont \S 2. A reductive p-adic context}

\bigskip
In 2.1 we list conventions and notations. In 2.2 we recall a problem and some results. In 2.3 and 2.4 we include solutions of this problem for certain situations involving reductive group schemes that are not semisimple. 

Let $q:=p^r$, where $r\in\dbN$. Let $k:=\dbF_{q}$. Always $n\in\dbN$. Let $W_n(k):=W(k)/p^nW(k)$. 

\bigskip\noindent
{\bf 2.1. Notations and conventions.} We abbreviate absolutely simple as a.s. and simply connected as s.c. Let $R$ and $F$ be as in the beginning of \S 1. We say $F^{\ad}$ is simple (resp. is a.s.) if (resp. if each geometric fibre of) it has no proper, normal subgroup of positive relative dimension. For $S$ a closed subgroup of $F$, let $\Lie(S)$ be its $R$-Lie algebra. If $R$ is an $\dbF_p$-algebra, let $\Lie_{\dbF_p}(S)$ be $\Lie(S)$ but viewed either as an abstract group or as an $\dbF_p$-Lie algebra. As sets, we identify $\Lie(S)=\Ker(S(R[x]/x^2)\to S(R))$, where the $R$-homomorphism $R[x]/(x^2)\twoheadrightarrow R$ takes $x$ into $0$. If $F_1\to F$ is an \'etale isogeny, then we identity canonically $\Lie(F_1)=\Lie(F)$. If $R_1\hookrightarrow R$ is a finite and flat $\dbZ$-monomorphism, then $\Res_{R/R_1} S$ is the group scheme over $R_1$ obtained from $S$ through the Weil restriction of scalars (see [7, 1.5] and [3, 7.6]). It is defined by the functorial identity $\Hom_{\Spec(R)}(Y\times_{\Spec(R_1)} \Spec(R),S)=\Hom_{\Spec(R_1)}(Y,\Res_{R/R_1} S)$, where $Y$ is an arbitrary $R_1$-scheme. We get a canonical identification $\Lie_{R_1}(S)=\Lie(\Res_{R/R_1} S)$, where $\Lie_{R_1}(S)$ is $\Lie(S)$ but viewed just as an $R_1$-Lie algebra. All modules over subgroups of $F(R)$ are left modules.

Let $F$ be semisimple. Let $o(F)$ be the order of $Z(F)$ as a finite, flat, group scheme. So $o(SL_2)=2$ and $c(F)=o(F^{\sc})/o(F)$. See [6] for the classification of connected Dynkin diagrams. We say $F$ is of isotypic $DT\in\{A_n,B_n,C_n|n\in\dbN\}\cup\{D_n|n\in\dbN\setminus\{1,2\}\}\cup\{E_6,E_7,E_8,F_4,G_2\}$ Dynkin type if the Dynkin diagram of every simple factor of a geometric fibre of $F^{\ad}$ is $DT$; if $F^{\ad}$ is a.s. we drop the word isotypic. We use the standard notations for classical reductive group schemes over $k$, $W(k)$, $\dbR$ or $\dbC$ (see [2]). So $PGL_n=GL_n^{\ad}=SL_n^{\ad}$, $PGSp_{2n}=Sp_{2n}^{\ad}=GSp_{2n}^{\ad}$, $PGU_n=SU_n^{\ad}$, etc.  Let 
$$F(R)^\prime:=\im(F^{\sc}(R)\to F(R)).$$
\noindent
{\bf 2.2. A review.} Let $G$ be a reductive group scheme over $W(k)$. As $Z(G)$ is a group scheme of multiplicative type (cf. [10, Vol. III, 4.1.7 of p. 173]), the group $Z(G)(k)$ has order prime to $p$. Let $\dbA^1$ be the affine line over $B(k):=W(k)[{1\over p}]$. The group $G(B(k))$ is endowed with the coarsest topology making all maps $G(B(k))\to\dbA^1(B(k))=B(k)$ induced by morphisms $G_{B(k)}\to\dbA^1$ to be continuous. We identify naturally $\Ker(G(W_{n+1}(k))\to G(W_n(k)))$ with $\Lie_{\dbF_p}(G_k)$. Let $K$ be a closed subgroup of $G(W(k))$. So $K$ is compact.

\medskip\noindent
{\bf Problem.} {\it Find practical conditions on $G$, $K$, and $p$ that imply $K=G(W(k))$.}

\medskip
This Problem was first considered for $SL_n$ and $Sp_{2n}$ groups over $\dbZ_p$ by Serre (see [28, IV] and [33, p. 52]). For $G$ semisimple, it is in essence solved in [37] (see 2.2.5). In 2.3 and 2.4 we present refinements of 2.2.5 for $G$ non-semisimple. We now recall some basic results. 
\medskip\noindent
{\bf 2.2.1. Lemma.} {\it {\bf (a)} We assume that either $p\Ge 3$ and $K$ surjects onto $G(W_2(k))$ or $p=2$ and $K$ surjects onto $G(W_3(k))$. Then $K=G(W(k))$.

\smallskip
{\bf (b)} We assume $p=2$. Let $x$, $y\in\Ker(G(W(k))\to G(k))$. Let $\bar x$, $\bar y\in\Ker(G(W_2(k))\to G(k))=\Lie_{\dbF_2}(G_k)$ be the reductions mod $4$ of $x$ and $y$. Then the image of the commutator $xyx^{-1}y^{-1}$ in $\Ker(G(W_3(k))\to G(W_2(k)))=\Lie_{\dbF_2}(G_k)$ is the Lie bracket $[\bar x,\bar y]$.}

\medskip
\proof
See [37, 4.1.2] for (a). Part (b) is just the computation [37, (20) of 4.7] performed for a reductive (instead of a semisimple) group scheme.\endproof

\medskip
See [37, 4.2 3)] for the following classical result. 

\medskip\noindent
{\bf 2.2.2. Proposition.} {\it The group $G^{\ad}$ (resp. $G^{\sc}$) is a product of Weil restrictions of a.s. adjoint groups (resp. of simply connected groups having a.s. adjoints) over Witt rings of finite field extensions of $k$.}
\medskip

Next we recall two well known results.

\medskip\noindent
{\bf 2.2.3. Proposition.} {\it We assume $G$ is semisimple. Let $S$ be a subgroup of $G(k)$ whose image in $G^{\ad}(k)$ contains $G^{\ad}(k)^\prime$. Then $G(k)^\prime$ is a subgroup of $S$.}

\medskip
\proof
We can assume that $S$ surjects onto $G^{\ad}(k)^\prime$. It is known that $G^{\sc}(k)$ is generated by elements of order $p$, cf. 2.2.2 and [15, 2.2.6 (f)]. So $G^{\ad}(k)^\prime$ is also generated by elements of order $p$. The kernel of the epimorphism $S\twoheadrightarrow G^{\ad}(k)^\prime$ is a subgroup of $Z(G)(k)$ and so of order prime to $p$. So any element of $G^{\ad}(k)^\prime$ of order $p$ is the image of an element of $S$ of order $p$. So by replacing $S$ by its subgroup generated by elements of order $p$, we can assume that $S$ is generated by elements of order $p$. So as the group $G(k)/G(k)^\prime$ has order prime to $p$, we have $S\leqslant G(k)^\prime$. So the inverse image of $S$ into $G^{\sc}(k)$ surjects onto $S$ and so also onto $G^{\ad}(k)^\prime$. 

So it suffices to prove the Lemma under the extra assumption $G=G^{\sc}$. We have $G(k)/Z(G)(k)=S/S\cap Z(G)(k)\times Z(G)(k)/S\cap Z(G)(k)$. So as $G(k)$ is generated by elements of order $p$ and as $Z(G)(k)$ has order prime to $p$, the group $Z(G)(k)/S\cap Z(G)(k)$ is trivial. So $Z(G)(k)\vartriangleleft S$ and so $S=G(k)$.\endproof 

\medskip\noindent
{\bf 2.2.4. Theorem.} {\it If $q=3$ we assume $G^{\ad}$ has no simple factor that is a $PGL_2$ group, and if $q=2$ we assume $G^{\ad}$ has no simple factor that is a $PGL_2$, $PGSp_4$, $PGU_3$ or a split group of $G_2$ Dynkin type. Then the group $G^{\ad}(k)^\prime$ is a product of non-abelian simple groups generated by elements of order $p$; moreover, if $G^{\ad}$ is simple, then $G^{\ad}(k)^\prime$ is non-abelian simple. So all factors of the composition series of $G^{\der}(k)$ are either non-abelian simple or cyclic of order prime to $p$.}

\medskip
\proof
We can assume $G^{\ad}$ is a.s., cf. 2.2.2. So the first part follows from [15, 2.2.1, 2.2.6 (f) and 2.2.7 (a)]. We have short exact sequences $0\to Z(G^{\sc}_k)(k)\to G^{\sc}(k)\to G^{\ad}(k)^\prime\to 0$ and $0\to G^{\der}(k)^\prime\to G^{\der}(k)\to H^1(k,\Ker(G^{\sc}_k\to G^{\der}_k))\to 0$ and the finite abelian groups $Z(G^{\sc})(k)$ and $H^1(k,\Ker(G^{\sc}_k\to G^{\der}_k))$ are of order prime to $p$. From this and the first part we get the second part.\endproof

\medskip
The next Theorem is a slightly weaker form of the combined results [37, 1.3, 4.1.1, 4.5] to be used often in what follows. 

\medskip\noindent
{\bf 2.2.5. Theorem.} {\it We assume that $G$ is semisimple. Then $G(W(k))$ is the unique closed subgroup of $G(W(k))$ surjecting onto $G(k)$ iff the following two statements hold:

\medskip
{\bf (i)} we have $g.c.d.(p,c(G))=1$, and 

\smallskip
{\bf (ii)} if $q\in\{3,4\}$, then $G^{\ad}$ has no simple factor that is a $PGL_2$ group, and if $q=2$, then $G^{\ad}$ has no simple factor that is a $PGL_2$, $PGL_3$, $PGU_3$, $PGU_4$, $\Res_{W(\dbF_4)/\dbZ_2} PGL_2$ or a split group of $G_2$ Dynkin type.}

\medskip
\proof
The ``if" part is just a weaker form of [37, Th. 1.3]. 
We now check the ``only if" part. So $G(W(k))$ has no proper, closed subgroup surjecting onto $G(k)$. This implies that the epimorphism $G(W_2(k))\twoheadrightarrow G(k)$ has no right inverse. If (i) does not hold, then there is an isogeny cover of $G$ of degree $p$ and so from [37, Ex. 4.1.1] we get that $G(W(k))$ has proper, closed subgroups surjecting onto $G(k)$. So (i) holds. We show that the assumption that (ii) does not hold leads to a contradiction. 

So $q\in\{2,3,4\}$ and $G^{\ad}$ has a simple factor $G_0$ among those listed in (ii). Let $F$ be a semisimple, normal, closed subgroup of $G$ such that $(G/F)^{\ad}$ is naturally isomorphic to $G_0$. All torsors of $F_k$ and so also of $F$ are trivial (see [32, p. 132]). So we have a short exact sequence $0\to F(W(k))\to G(W(k))\to (G/F)(W(k))\to 0$. So the inverse image into $G(W(k))$ of any proper, closed subgroup of  $(G/F)(W(k))$ surjecting onto $(G/F)(k)$, is a proper, closed subgroup of $G(W(k))$ surjecting onto $G(k)$. So in order to reach a contradiction, we can assume $F$ is the trivial subgroup of $G$. So $G_0=G^{\ad}$. So from this and (i) we get that either $q=4$ and $G$ is an $SL_2$ group, or $q=3$ and $G$ is a $PGL_2$ or an $SL_2$ group, or $q=2$ and $G$ is either a split group of $G_2$ Lie type or an $SL_2$, $SL_3$, $PGL_3$, $SU_3$, $PGU_3$, $SU_4$ or a $\Res_{W(\dbF_4)/\dbZ_2} SL_2$ group. 

So as the epimorphism $G(W_2(k))\twoheadrightarrow G(k)$ has no right inverse, from [37, Th. 4.5] we get that either $q=4$ and $G$ is an $SL_2$ group or $q=2$ and $G$ is an $SL_2$, $SU_4$ or a $\Res_{W(\dbF_4)} SL_2$ group. It is well known that $G$ can not be an $SL_2$ group over $\dbZ_2$ and so from now on we will exclude this case. In the remaining three cases with $q\in\{2,4\}$, the group $G^{\ad}(k)^\prime=G^{\ad}(k)=G(k)$ is either $SL_2(\dbF_4)$ or $SU_4(\dbF_2)$ and so it is non-abelian simple (cf. 2.2.4). Let $\Sigma$ be a subgroup of $G^{\ad}(W_2(k))$ surjecting isomorphically onto $G^{\ad}(k)$, cf. [37, Th. 4.5]. Let $K$ be the inverse image of $\Sigma$ via the homomorphism $G(W(k))\to G^{\ad}(W_2(k))$. It is a proper, closed subgroup of $G(W(k))$ containing $\Ker(G(W(k))\to G(W_2(k)))$. As $\Sigma$ is a non-abelian simple group and as $G^{\ad}(W_2(k))/\im(G(W_2(k))\to G^{\ad}(W_2(k)))$ is an abelian $2$-group, we easily get that $\Sigma$ is a subgroup of $\im(G(W_2(k))\to G^{\ad}(W_2(k)))$. So $K$ surjects onto $G(k)$. Contradiction. So (ii) holds.\endproof

\bigskip\noindent
{\bf 2.3. Proposition.} 
{\it Let $F$ be a semisimple, normal, closed subgroup of $G^{\der}$. Let $K$ be a closed subgroup of $G(W(k))$. Let $K_1:=\im(K\to G(k))$. Let $K^{\d}:=K\cap F(W(k))$. Let $K_1^{\d}:=\im(K^{\d}\to F(k))$. We assume that $g.c.d.(p,c(F))=1$ and that $F(k)$ (resp. and that $F(k)^\prime$) is a subgroup of $K_1$. If $q\in\{3,4\}$ we also assume that $F^{\ad}$ has no simple factor that is a $PGL_2$ group, and if $q=2$ we also assume that $F^{\ad}$ has no simple factor that is a $PGL_2$, $PGL_3$, $PGSp_4$, $PGU_3$, $PGU_4$, $\Res_{W(\dbF_4)/\dbZ_2} PGL_2$ or a split group of $G_2$ Dynkin type. Then we have $K^{\d}=F(W(k))$ (resp. we have $F(k)^\prime\vartriangleleft K_1^{\d}$ and $\Ker(F(W(k))\to F(k))\vartriangleleft K^{\d}$).}

\medskip
\proof
Let $g\in K_1\cap F(k)$ be an arbitrary element. As we have $F(k)^\prime\vartriangleleft K_1\cap F(k)\vartriangleleft F(k)$, the quotient group $F(k)/K_1\cap F(k)$ is abelian. Moreover, all factors of the composition series of $K_1\cap F(k)$ are either non-abelian simple or cyclic of order prime to $p$ (cf. 2.2.4). So we can write $g$ as a product $g_1^{\prime}g_2^{\prime}...g_u^{\prime}$, where $u\in\dbN$ and each $g_i^\prime$ is either of the form $g_i^p$ or of the form $[g_{i_1},g_{i_2}]$, for some $g_i$, $g_{i_1}$, $g_{i_2}\in K_1\cap F(k)$ ($i\in\{1,...,u\}$). 

We check by induction on $s\in\dbN$ that for any $g\in K_1\cap F(k)$ there is $h_s(g)\in K\cap\Ker(G(W(k))\to (G/F)(W_s(k)))$ whose reduction mod $p$ is the element $g$ of $G(k)$. The existence of $h_1(g)$ follows from the definition of $K_1$. The passage from $s$ to $s+1$ goes as follows. If $g_i^\prime=g_i^p$, let $h_i^{\prime}(g):=h_s(g_i)^p$. If $g_i=[g_{i_1},g_{i_2}]$, let $h_i^\prime(g):=[h_s(g_{i_1}),h_s(g_{i_2})]$. As $\Ker((G/F)(W_{s+1}(k))\to (G/F)(W_s(k)))$ is an abelian $p$-group, for $i\in\{1,...,u\}$ we have $h_i^\prime(g)\in K\cap\Ker(G(W(k))\to (G/F)(W_{s+1}(k)))$. So we can take $h_{s+1}(g)$ to be the product $h_1^{\prime}(g)h_2^{\prime}(g)...h_u^{\prime}(g)$. This ends the induction.

As $K$ is compact, there is a subsequence of the sequence $(h_s(g))_{s\in\dbN}$ that converges to an element $h(g)\in K^{\d}$. The reduction mod $p$ of $h(g)$ is $g$. 

So if $F(k)\vartriangleleft K_1$, then $g$ is an arbitrary element of $F(k)$ and so $K^{\d}$ surjects onto $F(k)$. Moreover from 2.2.5 we get that $K^{\d}=F(W(k))$. 

We now consider the case when $F(k)^\prime\vartriangleleft K_1$. By taking $g$ to be an arbitrary element of $F(k)^\prime$, from the existence of $h(g)\in K^d$ we get $F(k)^\prime\vartriangleleft K_1^{\d}$. 

Let $K^{\s}$ be the inverse image of $K^{\d}$ into $F^{\sc}(W(k))$. As $g.c.d.(p,c(F))=1$, the isogeny $F^{\sc}\to F$ is \'etale. So $\Ker(F(W(k))\to F(k))=\Ker(F^{\sc}(W(k))\to F^{\sc}(k))$ and the group $K^{\d}/\im(K^{\s}\to K^{\d})$ is finite of order prime to $p$. So the image of $\im(K^{\s}\to F(k))$ in $F(k)^\prime$ is a normal subgroup of index prime to $p$ and so it is $F(k)^\prime$, cf. 2.2.4. So $\im(K^{\s}\to F(k))=F^{\sc}(k)$, cf. 2.2.3. So $K^{\s}=F^{\sc}(W(k))$, cf 2.2.5. So $\Ker(F(W(k))\to F(k))\vartriangleleft \im(K^{\s}\to K^{\d})\vartriangleleft K^{\d}$.\endproof

\bigskip\noindent 
{\bf 2.4. A setting.} Let $\tilde G:=G/Z^0(G)$ (see beginning of \S1 for $Z^0(G)$). The hypothesis $g.c.d.(p,c(G^{\der}))=1$ of 2.3 (applied with $F=G^{\der}$) is often too restrictive. So in this section we will present a new approach that on one side allows us to treat cases with $p$ dividing $c(G^{\der})$ and on another side gives us criteria of when $K$ surjects onto $\tilde G(W(k))$. These last criteria are especially of interest when $p$ divides $c(\tilde G)$ (i.e. when we can not appeal to 2.2.5). Warning: until \S3 we will use and refer to the conditions (i) to (v) below.

\medskip
{\bf (i)} The adjoint group $\tilde G^{\ad}$ is simple and we have $g.c.d.(p,o(\tilde G))=1$.  

\smallskip
{\bf (ii)} The factors of the composition series of $\tilde G(k)$ are either cyclic of order prime to the order of $G^{\ab}(W_2(k))$ or non-abelian simple groups.

\smallskip
{\bf (iii)} We can write the natural isogeny $Z^0(G)\to G^{\ab}$ as a product $\prod_{i\in I} h_i:T_i\to T_i^\prime$ of isogenies in such a way that there is a non-empty subset $I_p$ of $I$ with the property that for any $i\in I\setminus I_p$ (resp. for any $i\in I_p$), the isogeny $h_i:T_i\to T_i^\prime$ has degree prime to $p$ (resp. factors through the $p$-th power endomorphism of $T_i$).

\smallskip
{\bf (iv)} We have $\sum_{i\in I_p} \dim_k(T_{ik})=\dim_k(\Lie(\tilde G_k)/[\Lie(\tilde G_k),\Lie(\tilde G_k)])$.

\medskip
Condition (i) in essence says that $\tilde G$ ``behaves" like a simple, adjoint group. Condition (ii) is almost always satisfied, cf. 2.2.4. Conditions (iii) and (iv) require that the tori $Z^0(G)$ and $G^{\ab}$ have big enough ranks and the isogeny $Z^0(G)\to G^{\ab}$ has degree divisible by a sufficiently big power of $p$. So roughly speaking, in order that (iii) and (iv) hold, we will use reductive group schemes $G$ for which the isogeny $Z^0(G)\to G^{\ab}$ ``behaves" like the analogue isogeny for the case of a (Weil restriction of a) $GL_{pn}$ group.

Let $k_1$ be the finite field extension of $k$ such that $\tilde G^{\ad}$ is the $\Res_{W(k_1)/W(k)}$ of an a.s. adjoint group $\tilde G_1^{\ad}$ over $W(k_1)$, cf. 2.2.2. The next condition will be needed only if $p=2$.

\medskip
{\bf (v)} If $p=2$, then at least one of the following two assumptions holds:

\medskip
\item{{\bf (va)}} the image of $\prod_{i\in I_p} T_i(W(k))$ in $G^{\ab}(W_3(k))$ has odd order (for instance, this holds if $q=2$ and $\prod_{i\in I_p} T_i$ is a split torus);

\smallskip
\item{{\bf (vb)}} the field $k_1$ has at least four elements and $\tilde G^{\ad}_1$ is either split or a $PGU_{2n+2}$ group.

\medskip\noindent
{\bf 2.4.1. Theorem.} {\it Let $G$ be a reductive group scheme over $W(k)$ such that all the above conditions (i) to (v) hold. Let $K$ be a closed subgroup of $G(W(k))$ surjecting onto $\tilde G(k)$ as well as onto $G^{\ab}(W(k))$. Depending on the parity of $p$ we have:

\medskip
{\bf (a)} Let $p>2$. Then $K$ surjects onto $\tilde G(W(k))$. If moreover we have $I=I_p$, then the images of $K$ and $G(W(k))$ in $\tilde G(W_2(k))\times G^{\ab}(W_2(k))$  are the same and isomorphic to $\tilde G(W_2(k))\times G^{\ab}(k)$.

\smallskip
{\bf (b)} Let $p=2$. Let $\tilde K_3^2$ be the intersection of $\im(K\to \tilde G(W_3(k)))$ with $\Ker(\tilde G(W_3(k))\to\tilde G(W_2(k)))=\Lie_{\dbF_p}(\tilde G_k)=\Lie_{\dbF_p}(\tilde G_{1k_1}^{\ad})$. If we work under (vb), then we also assume that $\tilde K_3^2$ is a $k_1$-vector subspace of $\Lie(\tilde G_{1k_1}^{\ad})$. Then $K$ surjects onto $\tilde G(W(k))$. If moreover we have $I=I_2$ and we work under (va) (resp. and we work under (vb)), then the images of $K$ and $G(W(k))$ in $\tilde G(W_s(k))\times G^{\ab}(W_s(k))$, with $s=3$ (resp. with $s=2$), are the same and isomorphic to $\tilde G(W_s(k))\times G^{\ab}(k)$.} 

\medskip
\proof
We can assume $I=I_p$ (otherwise we replace $G$ by $G/\prod_{i\in I\setminus I_p} T_i$). From (iii) we get that $p|o(G^{\der})$. So $\tilde G$ is not of isotypic $E_8$, $F_4$ or $G_2$ Dynkin type and the right hand side of (iv) is at least 1. So the torus $G^{\ab}$ is nontrivial. Let $K_{22}$ (resp. $K_{21}$) be the image of $K$ in $\tilde G(W_2(k))\times G^{\ab}(W_2(k))$
(resp. in $\tilde G(W_2(k))\times G^{\ab}(k)$). From (ii) and our surjectivity hypotheses we get that $K$ surjects onto $\tilde G(k)\times G^{\ab}(k)$. So
we have a short exact sequence
$$0\to K_{21}^{11}\to K_{21}\to \tilde G(k)\times G^{\ab}(k)\to 0,$$
where $K_{21}^{11}:=K_{21}\cap (\Ker(\tilde G(W_2(k))\times G^{\ab}(k)\to \tilde G(k)\times G^{\ab}(k)))=K_{21}\cap\Lie_{\dbF_p}(\tilde G_k)$. Let
$K_{22}^{21}:=K_{22}\cap (\Ker(\tilde G(W_2(k))\times G^{\ab}(W_2(k))\to\tilde G(W_2(k))\times G^{\ab}(k)))=K_{22}\cap\Lie_{\dbF_2}(G^{\ab}_k)$. 
We have a second short exact sequence 
$$0\to K_{22}^{21}\to K_{22}\to K_{21}\to 0.$$ 
\indent
The isogeny $Z^0(G)\to G^{\ab}$ factors through the $p$-th power endomorphism of $Z^0(G)$ (cf. (iii)) and moreover the group $Z^0(G)(k)$ has order prime to $p$. So the homomorphism $\Ker(Z^0(G)(W_2(k))\to G^{\ab}(k))\to G^{\ab}(W_2(k))$ is trivial. So as we have a short exact sequence $0\to Z^0(G)(W_2(k))\to G(W_2(k))\to\tilde G(W_2(k))\to 0$, the group $K_{22}^{21}$ is trivial.
So we can identify $K_{22}=K_{21}$ and so we will also view $K_{21}^{11}$ as a subgroup of $\tilde G(W_2(k))\times G^{\ab}(W_2(k))$. Each homomorphism $\tilde G(k)\to G^{\ab}(W_2(k))/\im(K_{21}^{11}\to G^{\ab}(W_2(k)))$ is trivial, cf. (ii). As $K$ surjects onto $G^{\ab}(W_2(k))$, from the last two sentences we get that $K_{21}^{11}$ surjects onto $\Ker(G^{\ab}(W_2(k))\to G^{\ab}(k))$ and so is a nontrivial abelian group. As $K$ surjects onto $\tilde G(k)$, under the adjoint representation $\tilde G(k)\to GL(\Lie_{\dbF_p}(\tilde G_k))$ both $K_{21}^{11}$ and $\Ker(K_{21}^{11}\to \Ker(G^{\ab}(W_2(k))\to G^{\ab}(k)))$ are $\tilde G(k)$-modules. Let $\grL:=[\Lie_{\dbF_p}(\tilde G_k),\Lie_{\dbF_p}(\tilde G_k)]$.

We first consider the case when either $p\Ge 3$ or $p=2$ and $\tilde G$ is not of isotypic $B_{n}$ or $C_{n}$ Dynkin type with $n\Ge 2$. The $\tilde G(k)$-module $\grL$ is simple and moreover $\tilde G(k)$ acts trivially on $\Lie_{\dbF_p}(\tilde G_k)/\grL$, cf. [37, 3.7.1 and 3.10 5)] applied to $\Res_{k/\dbF_p} \tilde G_k$. From this and [37, 3.10 1)] we get that the only simple $\tilde G(k)$-submodule of $\Lie_{\dbF_p}(\tilde G_k)$ is $\grL$. So $\grL\vartriangleleft K_{21}^{11}$. The image of $\grL$ in $\Ker(G^{\ab}(W_2(k))\to G^{\ab}(k))$ is trivial. So as $K_{21}^{11}$ surjects onto $\Ker(G^{\ab}(W_2(k))\to G^{\ab}(k))$, the group $K_{21}^{11}/\grL$ has at least as many elements as $\Lie_{\dbF_p}(G^{\ab}_k)$ and so (cf. (iv)) as $\Lie_{\dbF_p}(\tilde G_k)/\grL$. By reasons of orders we get $K_{21}^{11}=\Lie_{\dbF_p}(\tilde G_k)$. 

We now consider the case when $p=2$ and $\tilde G$ is of isotypic $B_{n}$ or $C_{n}$ Dynkin type with $n\Ge 2$. The group $\tilde G_1^{\ad}$ is not a $PGSp_4$ group over $\dbZ_2$, cf. (ii) and [15, 2.2.7 (a)]. We have $o(G^{\sc})=2^{[k_1:k]}$ and so from (i) we get that $o(\tilde G)=1$. Thus $\tilde G=\tilde G^{\ad}$. So $\dim_k(\Lie(\tilde G_k)/[\Lie(\tilde G_k),\Lie(\tilde G_k)])=[k_1:k]\dim_{k_1}(\Lie(\tilde G_{1k_1}^{\ad})/[\Lie(\tilde G_{1k_1}^{\ad}),\Lie(\tilde G_{1k_1}^{\ad})])=[k_1:k]$, cf. [16, $(B_r)$ and $(C_r)$ of 0.13]. So as $\prod_{i\in I} \Ker(h_i)$ is a subgroup of $Z(G^{\der})$, from (iii) and (iv) we get that $o(G^{\der})$ is at least $2^{[k_1:k]}$. Thus $o(G^{\der})=o(G^{\sc})=2^{[k_1:k]}$ and so $G^{\der}$ is s.c. The finite group $G^{\der}(k)=\tilde G(k)$ is non-abelian simple (cf. 2.2.4). So as $G(k)=G^{\der}(k)\times G^{\ab}(k)$ and as $K$ surjects onto $G^{\der}(k)$, we have $G^{\der}(k)\vartriangleleft\im(K\to G(k))$. So $G^{\der}(W(k))\vartriangleleft K$, cf. 2.3 applied with $F=G^{\der}$. As $K$ surjects onto $G^{\ab}(W(k))$, we get $K=G(W(k))$. So we also have $K_{21}^{11}=\Lie_{\dbF_2}(\tilde G_k)$.

So regardless of who $p$ and $\tilde G$ are, we always have $K_{21}^{11}=\Lie_{\dbF_p}(\tilde G_k)$. So $K_{21}=\tilde G(W_2(k))\times G^{\ab}(k)$. But $K_{22}$ is isomorphic to $K_{21}$ and so also to $\tilde G(W_2(k))\times G^{\ab}(k)$. Applying this with $K$ being replaced by $G(W(k))$, we get that that the images of $K$ and $G(W(k))$ in $\tilde G(W_2(k))\times G^{\ab}(W_2(k))$ are the same. In particular $K$ surjects onto $\tilde G(W_2(k))$. So if $p\Ge 3$, then $K$ surjects onto $\tilde G(W(k))$ (cf. 2.2.1 (a)). This proves (a).

We now take $p=2$ and we prove the remaining part of (b). As $o(G^{\der})$ is even, the group $\tilde G^{\ad}$ is of isotypic $A_{2n+1}$, $B_n$, $C_n$, $D_n$ or $E_7$ Dynkin type. So as $K$ surjects onto $\tilde G(W_2(k))$, from 2.2.1 (b) we get that $\grL\leqslant \tilde K_3^2\leqslant \Lie_{\dbF_2}(\tilde G_k)$.

If (va) holds, then we proceed as in the mod $4$ context. Let $K_{33}$ (resp. $K_{31}$) be the image of $K$ in $\tilde G(W_3(k))\times G^{\ab}(W_3(k))$ (resp. in $\tilde G(W_3(k))\times G^{\ab}(k)$). Let $K_{33}^{22}:=\Ker(K_{33}\to\tilde G(W_2(k))\times G^{\ab}(W_2(k)))$. Let $\tilde g_3^2\in\tilde K_3^2$. Let $g_3\in K_{33}$ surjecting onto $\tilde g_3^2$. As $K_{22}\arrowsim\tilde G(W_2(k))\times G^{\ab}(k)$ and as the order of $G^{\ab}(k)$ is prime to $p$, the image of $g_{33}$ in $K_{22}$ is the identity element and so $g_{33}\in K_{33}^{22}$. So the natural homomorphism $K_{33}^{22}\to\tilde K_3^2$ is surjective. So based on (va), in the same way we got that we can identity $K_{22}=K_{21}$, we now get that we can identify $K_{33}^{22}=\tilde K_3^2$ and $K_{33}=K_{31}$. As $K$ surjects onto $G^{\ab}(W_3(k))$, the group $K_{33}^{22}$ surjects onto $\Ker(G^{\ab}(W_3(k))\to G^{\ab}(W_2(k)))=\Lie(G^{\ab}_k)$. If $\tilde G^{\ad}$ is not of isotypic $B_n$ or $C_n$ Dynkin type with $n\Ge 2$, then entirely as in the mod $4$ context for $K_{21}^{11}$ we argue that $K_{33}^{22}=\tilde K_3^2$ is $\Lie_{\dbF_2}(\tilde G_k)$. If $\tilde G^{\ad}$ is of isotypic $B_n$ or $C_n$ Dynkin type with $n\Ge 2$, then $K=G(W(k))$ (see above) and so we again have  $K_{33}^{22}=\tilde K_3^2=\Lie_{\dbF_2}(\tilde G_k)$. 

So $K_{33}=K_{31}=\tilde G(W_3(k))\times G^{\ab}(k)$. Applying this with $K$ being replaced by $G(W(k))$, we get that that the images of $K$ and $G(W(k))$ in $\tilde G(W_3(k))\times G^{\ab}(W_3(k))$ are the same and isomorphic to $\tilde G(W_3(k))\times G^{\ab}(k)$.

We now consider the case when (vb) holds. If $\tilde G$ is (resp. is not) of isotypic $D_{2n}$ Dynkin type, let $u:=2$ (resp. $u:=1$). We check that the assumptions of (vb) on $\tilde G_1^{\ad}$ imply the existence of $\dbG_m$ subgroups $T(j)$ of $\tilde G^{\ad}_1$, $j=\overline{1,u}$, such that the natural $k_1$-linear map 
$$l:\oplus_{j=1}^u\Lie(T(j)_{k_1})\to \Lie(\tilde G^{\ad}_{1k_1})/[\Lie(\tilde G^{\ad}_{1k_1}),\Lie(\tilde G^{\ad}_{1k_1})]$$ 
is an isomorphism. The codomain of $l$ has dimension $u$ (cf. [16, 0.13]) and the Lie algebra of any maximal torus of a Borel subgroup of $\tilde G^{\ad}_{1k_1}$ surjects onto it (cf. [37, 3.7 3)]). So the existence of $T(j)$'s when $\tilde G_1^{\ad}$ is a $PGU_{2n+2}$ group (resp. is split) is a consequence of the fact that the Lie algebra of any maximal split torus of the simply connected group cover $\tilde G_{1k_1}^{\sc}$ of $\tilde G_{1k_1}^{\ad}$ contains $\Lie(Z(\tilde G_{1k_1}^{\sc}))$ (resp. is obvious). 

As $k_1$ has at least four elements, any closed subgroup of $T(j)(W(k_1))$ surjecting onto $T(j)(W_2(k_1))$ has a nontrivial image in $\Ker(T(j)(W_3(k_1))\to T(j)(W_2(k_1)))$. From this and the facts that $\grL\leqslant\tilde K_3^2$ and that $\Lie_{\dbF_2}(T(j)_{k_1})=\Ker(T(j)(W_2(k_1))\to T(j)(k_1))\leqslant \im(K\to\tilde G^{\ad}(W_2(k)))$, we get that $\tilde K_3^2/\grL$ has non-zero elements of the form $l(x_j)$, where $x_j\in\Lie_{\dbF_2}(T(j)_{k_1})$. So the $k_1$-vector subspace $\tilde K_3^2/\grL$ of $\Lie(\tilde G_{1k_1})/[\Lie(\tilde G_{1k_1}),\Lie(\tilde G_{1k_1})]=\Lie(\tilde G^{\ad}_{1k_1})/[\Lie(\tilde G^{\ad}_{1k_1}),\Lie(\tilde G^{\ad}_{1k_1})]$,  contains $l(k_1x_j)$. So $\tilde K_3^2=\Lie_{\dbF_2}(\tilde G_k)=\Lie_{\dbF_2}(\tilde G_{1k_1}^{\ad})$. 

In both situations (va) and (vb) we got that $K$ surjects onto 
$\tilde G(W_3(k))$ and so onto $\tilde G(W(k))$, cf. 2.2.1 (a). So (b) holds.\endproof 

\medskip\noindent
{\bf 2.4.2. Variants.} {\bf (a)} We assume that conditions (i), (iii), (iv) and (v) hold. We also assume the following weaker variant of (ii): the factors of the composition series of $\tilde G(k)$ are either cyclic of order prime to $p$ or non-abelian simple groups. We also assume that $\Ker(G^{\ab}(W(k))\to G^{\ab}(k))\vartriangleleft \im(K\to G^{\ab}(W(k)))$. Then the proof of 2.4.1 applies to give us that $K$ surjects onto $\tilde G(W(k))$. 

\smallskip
{\bf (b)} We assume that either $p=2$ and $\tilde G$ is of isotypic $D_{n+3}$ Dynkin type or $p$ is arbitrary and $\tilde G$ is of isotypic $A_{p^2n-1}$ Dynkin type. If $q=2$ we also assume that $\tilde G^{\ad}$ is not a $PGU_4$ group. Then we have a variant of 2.4.1 in which we replace (i) by the condition:

\medskip
{\bf (i')} the group $\tilde G^{\ad}$ is simple and one of the following disjoint two conditions holds:

\medskip
\item{--} if $\tilde G$ is of isotypic  $A_{p^2n-1}$ Dynkin type, then we have an isogeny $\tilde G^{\prime}\to\tilde G$ of degree prime to $p$, where $\tilde G^{\prime}$ is the quotient of $\tilde G^{\sc}$ by the maximal finite, flat subgroup of $Z(\tilde G^{\sc})$ annihilated by $p$, 

\smallskip
\item{--} if $\tilde G$ is of isotypic $D_{n+3}$ Dynkin type, then $\tilde G$ is the $\Res_{W(k_1)/W(k)}$ of a $Spin_{n+3}^*/\mu_2$ group, where $*\in\{+,-\}$. 

\medskip
The only difference from the proof of 2.4.1 is in arguing that $K_{21}^{11}=\Lie_{\dbF_p}(\tilde G_k)$. Let $k_1$ as before (v). We have $\dim_k(\Lie(\tilde G_k)/[\Lie(\tilde G_k),\Lie(\tilde G_k)])=[k_1:k]$, cf. [16, $(A_r)$ and $(D_r)$ of 0.13]. The short exact sequence $0\to \Lie_{\dbF_p}(\tilde G_k^{\ad})\to \tilde G^{\ad}(W_2(k))\to\tilde G^{\ad}(k)\to 0$ does not have a section, cf. [37, Th. 4.5]. So the group $K_{21}^{11}$ is not included in $\Lie_{\dbF_p}(Z(\tilde G_k))$ and so $K_{21}^{11}$ has a nontrivial image in $\Lie_{\dbF_p}(\tilde G_k)/\Lie_{\dbF_p}(Z(\tilde G_k))$. But this last $\tilde G(k)$-module is canonically identified with a $\tilde G(k)$-submodule of $\Lie_{\dbF_p}(\tilde G_k^{\ad})$ containing $\im(\Lie_{\dbF_p}(\tilde G_k^{\sc})\to\Lie_{\dbF_p}(\tilde G_k^{\ad}))$ and so containing $[\Lie_{\dbF_p}(\tilde G_k^{\ad}),\Lie_{\dbF_p}(\tilde G_k^{\ad})]$ (cf. [37, 3.7 3)]). But $[\Lie_{\dbF_p}(\tilde G_k^{\ad}),\Lie_{\dbF_p}(\tilde G_k^{\ad})]$ is a simple $\tilde G(k)$-submodule of $\Lie_{\dbF_p}(\tilde G_k^{\ad})$ (cf. [37, 3.7.1]) and moreover it is the only simple $\tilde G(k)$-submodule of $\Lie_{\dbF_p}(\tilde G_k^{\ad})$ (cf. [37, 3.10 1) and 3.10 5)]). So the image of $K_{21}^{11}$ in $\Lie_{\dbF_p}(\tilde G_k^{\ad})$ contains $[\Lie_{\dbF_p}(\tilde G_k^{\ad}),\Lie_{\dbF_p}(\tilde G_k^{\ad})]$. 

If $\tilde G$ is either of isotypic $A_{p^2n-1}$ Dynkin type or of isotypic $D_{n+3}$ Dynkin type with $n$ even, then by reasons of dimensions we have $\Lie(\tilde G_k)/\Lie(Z(\tilde G_k))=[\Lie(\tilde G_k^{\ad}),\Lie(\tilde G_k^{\ad})]$ (cf. [16, $(A_r)$ and $(D_r)$ of 0.13]). So $K_{21}^{11}$ surjects onto $\Lie_{\dbF_p}(\tilde G_k)/\Lie_{\dbF_p}(Z(\tilde G_k))$. But $K_{21}^{11}$ surjects also onto the abelian group $\Lie_{\dbF_p}(G^{\ab}_k)$ having the same number of elements as $k_1$ (cf. (iv)) and so as $\Lie_{\dbF_p}(Z(\tilde G_k))$. Thus by reasons of orders of factors of composition series, we get that $K_{21}^{11}=\Lie_{\dbF_p}(\tilde G_k)$. 

Let now $\tilde G$ be of isotypic $D_{n+3}$ Dynkin type with $n$ odd; so we have $p=2$. As $[\Lie(\tilde G_k^{\ad}),\Lie(\tilde G_k^{\ad})]$ is a $k$-vector subspace of $\Lie(\tilde G_k)/\Lie(Z(\tilde G_k))$ of codimension $[k_1:k]$, the previous paragraph has to be slightly modified. The inverse image of $[\Lie(\tilde G_k^{\ad}),\Lie(\tilde G_k^{\ad})]$ in $\Lie(\tilde G_k)$ is $[\Lie(\tilde G_k),\Lie(\tilde G_k)]$ (to check this we can assume $k=k_1$ and this case follows from [16, $(D_r)$ of 0.13]). Moreover, $\grL$ is the only $\tilde G(k)$-submodule of itself surjecting onto $[\Lie_{\dbF_2}(\tilde G_k^{\ad}),\Lie_{\dbF_2}(\tilde G_k^{\ad})]$ (cf. [37, 3.10 3)]). Thus we have $\grL\vartriangleleft K_{21}^{11}$. So as in the proof of 2.4.1 (b) we get that $K_{21}^{11}=\Lie_{\dbF_2}(\tilde G_k)$. 

\medskip\noindent
{\bf 2.4.3. Examples.} 
{\bf (a)} Let $G$ be a $GL_{pn}$ group. So $\tilde G=\tilde G^{\ad}$ is a $PGL_{pn}$ group and so condition (i) holds. Moreover, we have $k_1=k$. If $q=2$ we assume that $n>1$. So we have the following weaker form of (ii): the factors of the composition series of $\tilde G(k)$ are either cyclic of order prime to $p$ or non-abelian simple groups (cf. 2.2.4). If $g.c.d.(n,q-1)=1$, then condition (ii) holds. The isogeny $Z^0(G)\to G^{\ab}$ can be identified with the $p$-power endomorphism of $\dbG_m$ and so by taking $I=I_p$ to have one element, we get that (iii) holds. We have $\sum_{i\in I_p} \dim_k(T_{ik})=1=\dim_k(\Lie(\tilde G_k)/[\Lie(\tilde G_k),\Lie(\tilde G_k)])$ (cf. [16, $(A_r)$ of 0.13] for the last equality). So (iv) holds. If $q=2$, then  (va) holds and if $p=2$ and $q>2$, then (vb) holds. So if $n$ is arbitrary (resp. if $g.c.d.(n,q-1)=1$), then the variant 2.4.2 (a) of 2.4.1 (resp. then 2.4.1) applies. 

\smallskip
{\bf (b)} Let $n_1$ be a divisor of $n\in\dbN$ that is prime to $p$. Let $m\in\dbN$ be at least $2$. Let $s\in\{1,...,m-1\}$. Let $G$ be a $GL_{p^mn}/\mu_{p^sn_1}$ group; so $g.c.d.(p,c(G^{\der}))=p>1$ and $\tilde G$ is a $PGL_{p^mn}$ group. If $n$ is arbitrary (resp. if $g.c.d.(n,q-1)=1$), then as in (a) we argue that the variant 2.4.2 (a) of 2.4.1 (resp. that 2.4.1) applies.

\smallskip
{\bf (c)} Let $p=2$. Let $G$ be a $GSpin_{2n+1}$ or a $GSp_{2n}$ group, with $n\Ge 2$; so $\tilde G$ is an adjoint group. If $q=2$ we assume that $G$ is not a $GSpin_5=GSp_4$ group. So $\tilde G(k)=\tilde G(k)^\prime$ is a non-abelian simple group, cf. 2.2.4. So condition (ii) holds. As in (a) we argue that conditions (i), (iii), (iv) and (v) hold. So 2.4.1 (b) applies.  

\smallskip
{\bf (d)} Let $*\in\{+,-\}$. Let $n\Ge 3$ be an odd integer. Let $p=2$. Let $G$ be a $GSO^*_{2n}$ group; so $g.c.d.(2,c(G^{\der}))=2>1$. The group $\tilde G$ is a simple, adjoint group of $D_n$ Dynkin type and so (i) holds. Also as $n\Ge 3$ the condition (ii) holds, cf. 2.2.4. The isogeny $Z^0(G)\to G^{\ab}$ can be identified with the $2$-power endomorphism of $\dbG_m$ and so by taking $I=I_2$ to have one element, we get that (iii) holds. We have $\sum_{i\in I_2} \dim_k(T_{ik})=1=\dim_k(\Lie(\tilde G_k)/[\Lie(\tilde G_k),\Lie(\tilde G_k)])$ (cf. [16, $(D_r)$ of 0.13] for the last equality; here is the place where we need $n$ to be odd). If $q=2$, then condition (va) holds. If $q>2$ and $*=+$, then condition (vb) holds. In particular, if $q=2$, then all conditions (i) to (v) hold and so 2.4.1 (b) applies. 

\smallskip
{\bf (e)} Let $*\in\{+,-\}$. Let $p=2$. Let $G$ be a $GSpin_{2n}^*$ group with $n\Ge 3$. So $\tilde G$ is an $SO^*_{2n}$ group and so condition (i') of 2.4.2 (b) holds. As in (d) we argue that conditions (ii) to (v) hold. So the variant 2.4.2 (b) of 2.4.1 applies. 

\medskip\noindent
{\bf 2.4.4. Remark.} If $1\Le \sum_{i\in I_p} \dim_{k}(T_{ik})<\dim_k(\Lie(\tilde G_k)/[\Lie(\tilde G_k),\Lie(\tilde G_k)])$, then one can use the proof of 2.4.1 to bound from above the order of $\Lie_{\dbF_p}(\tilde G_k)/K_{21}^{11}$.

\bigskip
\noindent
{\boldsectionfont \S3. Applications to the study of $K_p(A)$}

\bigskip
Let $E$, $i_E$, $A$, $d$, $H_1$, $H_A$, $p$, $T_p(A)$, $H_{A\dbZ_p}$, $\rho_{A,p}$, $G_p$, $E^{\text{conn}}$ and $K_p(A)$ be as in the beginning of \S 1. Let $\dbQ_{\mu_{p^{\infty}}}$ be the algebraic field extension of $\dbQ$ obtained by adjoining all $p$-power roots of $1$. Until end we study $K_p(A)$. In 3.1 and 3.2 we list the main assumptions and few simple properties. Warning: in 3.2 to 3.5 we assume that conditions 1.1 (a) to (c) hold and we apply 2.2 to 2.4. We end this Chapter with remarks (see 3.6). We recall that all simple factors of $H_{A\dbC}^{\ad}$ are of classical Lie type, cf. [26] or [8]. 

We have the following remarkable prediction (cf. [23] and [30]):

\bigskip\noindent
{\bf 3.1. Mumford--Tate conjecture for $(A,p)$.} {\it As subgroups of the group $GL(V_p(A))=GL(H_1\otimes_{\dbQ} \dbQ_p)$, we have $G_p=H_{A\dbQ_p}$.}
\medskip  

Until 3.6 we assume 3.1 holds (see [25, 5.14 and 5.15] and [38] for concrete cases).  

\medskip\noindent
{\bf 3.1.1. The $p$-torsion volume.} 
As 3.1 holds, $K_p(A)$ is an open subgroup of $H_{A\dbZ_p}(\dbZ_p)$. Let $\scrM(p)$ be the Haar measure on $H_{A\dbQ_p}(\dbQ_p)$ normalized by the fact that a compact, open subgroup of $H_{A\dbQ_p}(\dbQ_p)$ of maximum volume, has volume 1. We call the number 
$$\scrV_p(A):=\scrM(p)(K_p(A))\in (0,1]\cap\dbQ$$
as the $p$-torsion volume of (the isogeny class of) $A$. We think this is justified as for $n>>0$ we have $\scrV_p(A)=\scrM(p)(H_{A\dbZ_p}(\dbZ_p))[H_{A\dbZ_p}(\dbZ_p/p^n\dbZ_p):\im(K_p(A)\to H_{A\dbZ_p}(\dbZ_p/p^n\dbZ_p))]^{-1}$.

\medskip\noindent
{\bf 3.2. Simple properties.} 
Until 3.6 we assume that conditions 1.1 (a) to (c) hold, i.e. that 3.1 holds, that $H_{A\dbZ_p}$ is a reductive group scheme over $\dbZ_p$ and that by denoting 
$$T:=Z^0(H_{A\dbZ_p})\,\,\,\text{and}\,\,\,\tilde G:=H_{A\dbZ_p}/T$$ 
we have
$$\im(K_p(A)\to \tilde G(\dbF_p))=\tilde G(\dbF_p).\leqno (1)$$ 
We have $\scrM(p)(H_{A\dbZ_p}(\dbZ_p))=1$, cf. [35, 3.8.1 and 3.8.2]. So 
$$\scrV_p(A)^{-1}=[H_{A\dbZ_p}(\dbZ_p):K_p(A)]\in\dbN.\leqno (2)$$ 
As $H_A$ contains $Z(GL(H_1\otimes_{\dbZ}\dbQ))$, the torus $T$ contains the group scheme $Z(GL(T_p(A)))$ of homotheties of $T_p(A)=H_1\otimes_{\dbZ} \dbZ_p$. The group scheme $\tilde G$ is semisimple. 

We introduce three abstract groups:
$$K_p^T:=\Ker(K_p(A)\to \tilde G(\dbZ_p))=K_p(A)\cap T(\dbZ_p)\vartriangleleft T(\dbZ_p),$$
$$K_p^{\ab}:=\im(K_p(A)\to H_{A\dbZ_p}^{\ab}(\dbZ_p))\leqslant H^{\ab}_{A\dbZ_p}(\dbZ_p),\,\,\text{and}$$
$$K_p^{\d}:=\Ker(K_p(A)\to H^{\ab}_{A\dbZ_p}(\dbZ_p))=K_p(A)\cap H_{A\dbZ_p}^{\der}(\dbZ_p)\leqslant H_{A\dbZ_p}^{\der}(\dbZ_p).$$
All torsors of $T_{\dbF_p}$ and so also of $T$ are trivial (see [32, p. 132]). So the complex $0\to T(\dbZ_p)\to H_{A\dbZ_p}(\dbZ_p)\to \tilde G(\dbZ_p)\to 0$
is a short exact sequence. We get a first index formula 
$$\scrV_p(A)^{-1}=[T(\dbZ_p):K_p^T][\tilde G(\dbZ_p):\im(K_p(A)\to \tilde G(\dbZ_p))].\leqno (3)$$
\noindent
Similarly, the complex $0\to H_{A\dbZ_p}^{\der}(\dbZ_p)\to H_{A\dbZ_p}(\dbZ_p)\to H_{A\dbZ_p}^{\ab}(\dbZ_p)\to 0$ is a short exact sequence and so we have a second index formula
$$\scrV_p(A)^{-1}=[H_{A\dbZ_p}^{\der}(\dbZ_p):K_p^{\d}][H_{A\dbZ_p}^{\ab}(\dbZ_p):K_p^{\ab}].\leqno (4)$$
\noindent
{\bf 3.2.1. Lemma.} {\it We assume that $T$ is a $\dbG_m$ group, i.e. that $T$ acts as homotheties on $T_p(A)$. Then $T\cap H^{\der}_{A\dbZ_p}=\mu_2$. If moreover $E^{\text{conn}}$ is linearly disjoint from $\dbQ_{\mu_{p^{\infty}}}$, then $K_p^{\ab}=H^{\ab}_{A\dbZ_p}(\dbZ_p)$.}

\medskip
\proof
Let $h_A:\Res_{\dbC/\dbR} \dbG_m\to H_{A\dbR}$ be the homomorphism defining the Hodge $\dbZ$-structure on $H_1$. The image under $h_A$ of the compact subtorus of $\Res_{\dbC/\dbR} \dbG_m$ is contained in $H^{\der}_{A\dbR}$ and contains the $\mu_2$ subgroup of $Z(GL(H_1\otimes_{\dbZ} \dbR))$. So the $\mu_2$ subgroup of $T$ is also a subgroup of $H_{A\dbZ_p}^{\der}$. Let $p_A:T_p(A)\otimes_{\dbZ_p} T_p(A)\to\dbZ_p(1)=\text{proj.\, lim.}_{s\in\dbN} \mu_{p^s}(\overline{\dbQ_p})$ be the Weil pairing defined by a polarization of $A$. It is fixed by $H_{A\dbZ_p}^{\der}$. So $T\cap H^{\der}_{A\dbZ_p}$ is a subgroup of $\mu_2$. So $T\cap H^{\der}_{A\dbZ_p}=\mu_2$. So $H_{A\dbZ_p}^{\ab}(\dbZ_p)=(T/\mu_2)(\dbZ_p)$ acts faithfully on the $\dbZ_p$-span of $p_A$. As $E^{\text{conn}}$ is linearly disjoint from $\dbQ_{\mu_{p^{\infty}}}$, the Galois group $\Gal(\dbQ_{\mu_{p^{\infty}}}/\dbQ)$ is naturally a quotient of $\Gal(E^{\text{conn}})$. So the image of the $p$-adic cyclotomic character of $\Gal(E^{\text{conn}})$ can be identified with $\dbG_m(\dbZ_p)=H_{A\dbZ_p}^{\ab}(\dbZ_p)$ and so $K_p^{\ab}=H^{\ab}_{A\dbZ_p}(\dbZ_p)$.\endproof

\medskip
Let $Z:=\Ker(\tilde G^{\sc}\to H_{A\dbZ_p}^{\der})$. The order $o$ of the finite (abstract) group $Z(\dbF_p)$ is prime to $p$ (see first paragraph of 2.2). In 3.3 and 3.4 we consider two disjoint theories.

\bigskip\noindent
{\bf 3.3. The $g.c.d.(p,c(H_A^{\der}))=1$ theory.}
Until 3.4 we will consider only cases satisfying $g.c.d.(p,c(H_A^{\der}))=1$. So the finite, flat group scheme $Z$ of multiplicative type is also \'etale.

\medskip\noindent
{\bf 3.3.1. Theorem.} {\it We assume that $g.c.d.(p,c(H_A^{\der}))=1$ and that $T$ is a $\dbG_m$ group. If $p=3$ we also assume that $\tilde G^{\ad}$ has no simple factor that is a $PGL_2$ group and if $p=2$ we also assume that $\tilde G^{\ad}$ has no simple factor that is a $PGL_2$, $PGL_3$, $PGSp_4$, $PGU_3$, $PGU_4$ or a $\Res_{W(\dbF_4)/\dbZ_2} PGL_2$ group. We have the following four properties.

\medskip
 {\bf (a)} If $g.c.d.(p,c(\tilde G))=1$, then the following three statements are equivalent:

\medskip
{\bf (i)} $K_p(A)$ surjects onto $H_{A\dbZ_p}(\dbF_p)$;

\smallskip
{\bf (ii)} $K_p^T$ surjects onto $T(k)$;

\smallskip
{\bf (iii)} $\scrV_p(A)=p^{-m(A,p)}$, with $m(A,p)\in\dbN\cup\{0\}$.
\medskip

We now also assume that $E^{\text{conn}}$ is linearly disjoint from $\dbQ_{\mu_{p^{\infty}}}$. 

\smallskip
{\bf (b)} If $g.c.d.(p,c(\tilde G))=1$ and if for $p=2$ all simple factor of $\tilde G^{\ad}$ are of some isotypic $A_n$ Dynkin type, then $\scrV_p(A)^{-1}$ divides $2(p-1)$.

\smallskip
{\bf (c)} Always $\scrV_p(A)^{-1}$ is a divisor of $o$. So if $H_{A\dbZ_p}^{\der}$ is s.c., then $\scrV_p(A)=1$. 

\smallskip
{\bf (d)} If $K_p(A)$ surjects onto $H_{A\dbZ_p}(\dbF_p)$, then $\scrV_p(A)=1$.}

\medskip
\proof 
We prove (a). We have $\im(K_p(A)\to \tilde G(\dbZ_p))=\tilde G(\dbZ_p)$, cf. (1) and 2.2.5. So $\scrV_p(A)^{-1}=[T(\dbZ_p):K_p^T]$, cf. (3). So as $T(k)$ has order prime to $p$, (ii) and (iii) are equivalent. Obviously (i) implies (iii). If (ii) holds, then $\im(K_p(A)\to H_{A\dbZ_p}(\dbF_p))$ has $T(k)$ as a subgroup and surjects onto $\tilde G(\dbF_p)$; so (i) holds. So (a) holds. 

For the rest of the proof we will assume that $E^{\text{conn}}$ is linearly disjoint from $\dbQ_{\mu_{p^{\infty}}}$. So $K_p^{\ab}=H_{A\dbZ_p}^{\ab}(\dbZ_p)$ (cf. 3.2.1) and so $K_p(A)$ surjects onto $H_{A\dbZ_p}^{\ab}(\dbZ_p)$. We prove (b). Let $\tilde K$ be the subgroup of $\tilde G(\dbZ_p)\times H^{\ab}_{A\dbZ_p}(\dbZ_p)$ generated by $\im(K_p(A)\to \tilde G(\dbZ_p)\times H^{\ab}_{A\dbZ_p}(\dbZ_p))$ and by the elements of $\tilde G(\dbZ_p)\times H^{\ab}_{A\dbZ_p}(\dbZ_p)$ of order $p-1$ whose components in $\tilde G(\dbZ_p)$ are identity. We have $\Ker(H_{A\dbZ_p}\to \tilde G\times H^{\ab}_{A\dbZ_p})=T\cap H_{A\dbZ_p}^{\der}=\mu_2$. So to prove (b) it suffices to show that $\tilde K=\tilde G(\dbZ_p)\times H^{\ab}_{A\dbZ_p}(\dbZ_p)$. For this it is enough to show that $\tilde K$ surjects onto $\tilde G(\dbZ_p/p^s\dbZ_p)\times H^{\ab}_{A\dbZ_p}(\dbZ_p/p^s\dbZ_p)$, where $s=2$ if $p\Ge 3$ and $s=3$ if $p=2$ (cf. 2.2.1 (a)). As $\tilde K$ surjects onto $\tilde G(\dbF_p)\times H^{\ab}_{A\dbZ_p}(\dbF_p)$, $\tilde G(\dbZ_p)$ and $H^{\ab}_{A\dbZ_p}(\dbZ_p)$, we only have to show that $\Lie(\tilde G_{\dbF_p})$ has no $\dbF_p$-vector subspace of codimension 1 normalized by $\tilde G(\dbF_p)$. But this follows from [37, 3.7.1] (resp. [37, 3.10 2)]) applied to each factor $H$ of $\tilde G^{\sc}_{\dbF_p}$ having a simple adjoint not of (resp. of) isotypic $A_{pn-1}$ Dynkin type. So (b) holds.

We prove (c). We have $\scrV_p(A)^{-1}=[H_{A\dbZ_p}^{\der}(\dbZ_p):K_p^{\d}]$, cf. (4). As $K_p(A)$ surjects onto $\tilde G(\dbF_p)$, we also have $\tilde G^{\ad}(\dbF_p)^\prime\vartriangleleft \im(K_p(A)\to \tilde G^{\ad}(\dbF_p))$. So from 2.2.3 we get that $H_{A\dbZ_p}^{\der}(\dbF_p)^\prime\leqslant\im(K_p(A)\to H_{A\dbZ_p}^{\der}(\dbF_p))$. So as $g.c.d.(p,c(H_A^{\der}))=1$ we have $\Ker(H_{A\dbZ_p}^{\der}(\dbZ_p)\to H_{A\dbZ_p}^{\der}(\dbF_p))\vartriangleleft K_p^{\d}$, cf. 2.3 (applied with $F=G^{\der}$). So as we have a short exact sequence
$0\to H_{A\dbZ_p}^{\sc}(\dbF_p)/Z(\dbF_p)\to H_{A\dbZ_p}^{\der}(\dbF_p)\to H^1(\dbF_p,Z)\to 0$, the index $[H_{A\dbZ_p}^{\der}(\dbZ_p):K_p^{\d}]$ divides the order of $H^1(\dbF_p,Z)$ that is $o$. This proves (c). 

We prove (d). But $\scrV_p(A)^{-1}$ is prime to $p$ (cf. (c), as $g.c.d.(p,o)=1$) as well as a power of $p$ (cf. (a), as (iii) holds). So $\scrV_p(A)=1$ and so (d) holds.\endproof

\medskip\noindent
{\bf 3.3.2. Theorem.} {\it We assume that $p=3$ and that $g.c.d.(3,c(\tilde G))=1$. Then the index $[\tilde G(\dbZ_3):\im(K_3(A)\to \tilde G(\dbZ_3))]$ is $27^m$ for some $m\in\dbN\cup\{0\}$.}

\medskip
\proof 
We write $\tilde G^{\ad}=\tilde G_1^{\ad}\times \tilde G_2^{\ad}$, where $\tilde G_1^{\ad}$ is a product of $PGL_2$ groups and $\tilde G_2^{\ad}$ is a product of simple groups that are not $PGL_2$ groups. Let $\tilde G_2$ be the semisimple, normal, closed subgroup of $\tilde G$ naturally isogenous to $\tilde G_2^{\ad}$. As $g.c.d.(3,c(\tilde G))=1$ we have $g.c.d.(3,c(\tilde G_2))=1$. Also $\tilde G_2(\dbF_3)$ is a subgroup of $\im(K_3(A)\to \tilde G(\dbF_3))=\tilde G(\dbF_3)$ (cf. (1)). So $\tilde G_2(\dbZ_3)\vartriangleleft\im(K_3(A)\to \tilde G(\dbZ_3))$, cf. 2.3 (in the definition of $\tilde G_2$ we excluded the $PGL_2$ factors in order to be able to apply 2.3). 

We have a short exact sequence $0\to \tilde G_2(\dbZ_3)\to \tilde G(\dbZ_3)\to (\tilde G/\tilde G_2)(\dbZ_3)\to 0$. As the isogeny $\tilde G/\tilde G_2\to \tilde G_1^{\ad}$ is \'etale, we have $\Ker((\tilde G/\tilde G_2)(\dbZ_3)\to (\tilde G/\tilde G_2)(\dbF_3))=\Ker(\tilde G_1^{\ad}(\dbZ_3)\to \tilde G_1^{\ad}(\dbF_3))$. So to prove the Theorem it suffices to show that for any $s\in\dbN$, the subgroup $K^1_s:=\im(K_3(A)\to \tilde G_1^{\ad}(W_{s+1}(\dbF_3)))\cap\Ker(\tilde G_1^{\ad}(W_{s+1}(\dbF_3))\to \tilde G_1^{\ad}(W_s(\dbF_3)))$ of $\Lie(\tilde G^{\ad}_{1\dbF_3})=\Ker(\tilde G_1^{\ad}(W_{s+1}(\dbF_3))\to \tilde G_1^{\ad}(W_s(\dbF_3)))$ has index $27^{m(s)}$, where $m(s)\in\dbN\cup\{0\}$. The $\tilde G_1^{\ad}(\dbF_3)^\prime$-module $\Lie(\tilde G^{\ad}_{1\dbF_3})$ is a direct sum of simple $\tilde G_1^{\ad}(\dbF_3)^\prime$-modules that are Lie algebras of simple factors of $\tilde G_{1\dbF_3}^{\ad}$, cf. [37, 3.7.1] applied to the $A_1$ Dynkin type. As $K^1_s$ is a $\tilde G_1^{\ad}(\dbF_3)^\prime$-module (cf. (1)), it is isomorphic to a direct sum of such simple $\tilde G_1^{\ad}(\dbF_3)^\prime$-modules. So $[\Lie(\tilde G^{\ad}_{1\dbF_3}):K^1_s]$ is $27^{m(s)}$, where $m(s)\in\dbN\cup\{0\}$.\endproof

\bigskip\noindent
{\bf 3.4. The $g.c.d.(p,c(\tilde G))>1$ theory.} 
Until 3.6 we will consider only cases satisfying $g.c.d.(p,c(\tilde G))>1$. Sections 2.4.1 and 2.4.2 can be applied in many situations that either involve arbitrary primes $p$ and tori $H_A^{\ab}$ of arbitrary high dimension or are such that we can not apply 3.3.1 and 3.3.2 to them. To list few such applications we consider two cases.

\medskip\noindent
{\bf 3.4.1. Case 1.}
We assume $T$ is a $\dbG_m$ group. So $T\cap H_{A\dbZ_p}^{\der}=\mu_2$, cf. 3.2.1. So in order to be able to apply 2.4.1 and 2.4.2 we assume that $p=2$ and $H_A^{\ad}$ is a.s. of some Dynkin type $DT\in\{A_{2n+1}|n\in\dbN\}\cup\{B_n,C_n|n\in\dbN\setminus\{1,2\}\}\cup\{D_n|n\in\dbN\setminus\{1,2,3\}\}$. The $A_1$ and $B_2=C_2$ Dynkin types are excluded due to 2.4 (ii), cf. 2.2.4. If $DT=A_3$ we also assume $H_{A\dbZ_2}^{\der}$ is split and if $DT\in \{D_{2n}|n\in\dbN\setminus\{1\}\}$ we also assume $H_{A\dbZ_2}^{\der}$ is s.c. 

We now list few direct consequences of these assumptions. If $DT\in\{B_n,C_n|n\in\dbN\setminus\{1,2\}\}$, then $\tilde G$ is adjoint. If $DT\in\{A_{4n+1}|n\in\dbN\}\cup\{B_n,C_n|n\in\dbN\setminus\{1,2\}\}$, then $\Lie(\tilde G)=\Lie(\tilde G^{\ad})$. If $DT\notin\{A_{2n+1},D_{n}|n\in 1+2\dbN\}$, then $c(H_{A\dbZ_2}^{\der})$ is odd. If $c(H_{A\dbZ_2}^{\der})$ is odd, then $\Lie(H_{A\dbZ_2}^{\der})=\Lie(\tilde G^{\sc})$.

\medskip\noindent
{\bf 3.4.1.1. Theorem.}
{\it We also assume $E^{\text{conn}}$ is linearly disjoint from $\dbQ_{\mu_{2^{\infty}}}$. We have:

\medskip
{\bf (a)} The group $K_2(A)$ surjects onto $\tilde G(\dbZ_2)$ and the images of $K_2(A)$ and $H_{A\dbZ_2}(\dbZ_2)$ in $\tilde G(\dbZ_2/8\dbZ_2)\times H_{A\dbZ_2}^{\ab}(\dbZ_2/8\dbZ_2)$ are the same and isomorphic to $\tilde G(\dbZ_2/8\dbZ_2)$;

\smallskip
{\bf (b)} If $DT\in\{A_{4n-1},D_{2n+3}|n\in\dbN\}$ we assume $c(H_{A\dbZ_2}^{\der})$ is odd. Then $\scrV_2(A)=1$.}

\medskip
\proof
We have $K_2^{\ab}=H_{A\dbZ_2}^{\ab}(\dbZ_2)$, cf. 3.2.1. We check that the conditions needed to apply 2.4.1 (b) and 2.4.2 (b) hold for $H_{A\dbZ_2}$ and $K_2(A)$. If $DT\in\{A_{4n+3},D_{n+3}|n\in\dbN\}$ and $c(H_{A\dbZ_2}^{\der})$ is odd (resp. otherwise), then condition (i') of 2.4.2 (b) (resp. 2.4 (i)) is implied by our hypotheses on $DT$. Condition 2.4 (ii) holds, cf. 2.2.4 and the fact that $H_{A\dbZ_2}^{\ab}(\dbZ/4\dbZ)$ has order $2$. The isogeny $T\to H^{\ab}_{\dbZ_2}$ has degree $2$. So by taking $I=I_2$ to have one element we get that 2.4 (iii) holds. As $\dim_{\dbF_2}(\Lie(\tilde G_{\dbF_2})/[\Lie(\tilde G_{\dbF_2}),\Lie(\tilde G_{\dbF_2})])=1$ (cf. [16, $(A_r)$ to $(D_r)$ of 0.13]), 2.4 (iv) holds. As the torus $T$ is split, 2.4 (va) holds. So (a) follows from 2.4.2 (b) (resp. 2.4.1 (b)) applied to $H_{A\dbZ_2}$ and $K_2(A)$. 

We prove (b). As $K_2(A)$ surjects onto $H_{A\dbZ_2}(\dbF_2)=H_{A\dbZ_2}^{\der}(\dbF_2)=\tilde G(\dbF_2)$ and as $c(H_{A\dbZ_2}^{\der})$ is odd we have $K_2^{\d}=H_{A\dbZ_2}^{\der}(\dbZ_2)$, cf. 2.3. So (b) follows from (4).\endproof

\medskip\noindent
{\bf 3.4.2. Case 2.} We assume $T$ is not a $\dbG_m$ group, i.e. that the rank of $T$ is at least $2$. We do not know any literature studying when $K_p^{\ab}$ is $H_{A\dbZ_p}^{\ab}(\dbZ_p)$ for an absolutely simple abelian variety $A$.${}^1$ $\vfootnote{1} {The only exception is [33, Th. of p. 60] that implies that $\Ker(H_{A\dbZ_p}^{\ab}(\dbZ_p)\to H_{A\dbZ_p}^{\ab}(\dbF_p))$ is a subgroup of $K_p^{\ab}$ for $p>>0$.}$ So we do not know how to ``assure" in practice (by using some geometric or algebraic assumptions like the linear disjointness of $E^{\text{conn}}$ and $\dbQ_{\mu_{p^{\infty}}}$) that $K_p^{\ab}$ is $H_{A\dbZ_p}^{\ab}(\dbZ_p)$. So here we present only one situation to which 2.4.2 (a) applies. 

\medskip\noindent
{\bf 3.4.2.1. Theorem.} {\it We assume that the following three properties hold:
\medskip
{\bf (i)} the group $H_A^{\ad}$ is a.s. of $A_{pn-1}$ Dynkin type with $p+n>3$ and $H_{A\dbR}^{\ad}$ is a $PGU(a,pn-a)_{\dbR}$ group, for some $a\in\{1,...,pn-1\}$ prime to $p$;

\smallskip
{\bf (ii)} we have a product decomposition $H_{A\dbZ_p}=G\times_{\dbZ_p} T_1$,
with $G$ as a $GL_{pn}$ group and with $T_1$ as a $\dbG_m$ group (so $\tilde G=G^{\ad}=H_{A\dbZ_p}^{\ad}$);

\smallskip
{\bf (iii)} there is a prime $v$ of $E$ of index of ramification prime to $p$ and such that $A$ has semistable reduction with respect to it and $E^{\text{conn}}$ has a prime unramified over $v$. 

\medskip
Then $K_p(A)$ surjects onto $G(\dbZ_p)$ and so we have $\scrV_p(A)^{-1}=up^s$, where $u\in\dbN$ divides $(p-1)$ and $s\in\dbN\cup\{0\}$.} 

\medskip
\proof
We can assume $E=E^{\text{conn}}$. Let $E_v$ be the completion of $E$ with respect to $v$. We fix an embedding $i_{E_v}:E_v\hookrightarrow\dbC$ extending $i_E$ and we use it to identify naturally $E_v$ and $\overline{E_v}$ with subfields of $\dbC$. Let $\Gal(E_v):=\Gal(\overline{E_v}/E_v)$. Let $\rho_{\text{det}}:H_{A\dbQ_p}\twoheadrightarrow G^{\ab}_{\dbQ_p}$ be the natural epimorphism. See [13] and [14] for Fontaine comparison theory. 

When viewed as a 1-dimensional representation of $H_{A\dbQ_p}$, $\rho_{\text{det}}$ is a direct summand of the representation of $H_{A\dbQ_p}$ on the tensor algebra of $V_p(A)\oplus \Hom_{\dbQ_p}(V_p(A),\dbQ_p)$, cf. [40, 3.5]. Let $\rho_{\text{det},v}:\Gal(E_v)\to G^{\ab}_{\dbQ_p}(\dbQ_p)$ be the representation obtained by composing the natural homomorphism $\rho_{A,v}:\Gal(E_v)\to H_{A\dbQ_p}(\dbQ_p)$ with $\rho_{\text{det}}(\dbQ_p)$. As $A$ has semistable reduction with respect to $v$ (cf. (iii)), the representation of $\Gal(E_v)$ on $V_p(A)$ is semistable in the sense of [14] (cf. [36, Th. 0.2]). So as the category of semistable representations of $\Gal(E_v)$ is Tannakian (cf. [14, 5.17]), $\rho_{\text{det},v}$ is also a semistable representation of $\Gal(E_v)$. So there is $m\in\dbZ$ such that $\rho_{\text{det},v}$ can be identified with the $m$-th power of the $p$-adic cyclotomic character $\chi$ of $\Gal(E_v)$, cf. [14, 5.4.1]. One computes $m$ as follows. 

Let $\mu_A:\dbG_m\to GL(H_1\otimes_{\dbZ} \dbC)$ be as in \S 1. The embedding $i_{E_v}$ allows us to write $H_{A\dbC}$ as a product $G_{\dbC}\times_{\dbC} T_{1\dbC}$. We identify $G_{\dbC}$ with $GL(V)$, where $V$ is a complex vector space of dimension $pn$. Let $\mu_0:\dbG_m\to G_{\dbC}$ be the cocharacter that is the composite of the factorization of $\mu_A$ through $H_{A\dbC}$ with the projection epimorphism $H_{A\dbC}\twoheadrightarrow G_{\dbC}$. As any semistable representation of $\Gal(E_v)$ is also Hodge--Tate (see [14]), the composite of $\mu_0$ with the (determinant) epimorphism $G_{\dbC}\twoheadrightarrow G^{\ab}_{\dbC}$ can be identified with the $-m$-th power of the identical character of $\dbG_m$ (this fact is also a particular case of [41, Prop. of 1.3]). 

The cocharacter $\mu_{0}$ is a product of a cocharacter $\mu_{00}:\dbG_m\to Z(G_{\dbC})$ with a cocharacter $\mu_{01}:\dbG_m\to G_{\dbC}$ that acts trivially on a subspace $V_0$ of $V$ of dimension $c$ and via the identical character of $\dbG_m$ on $V/V_0$, where $c\in\{a,pn-a\}$ (cf. (i)). So $m$ is congruent mod $p$ to $-c$ and so to either $a$ or $-a$. So $g.c.d.(p,m)=1$, cf. (i). So as the index of $v$ is prime to $p$ (cf. (iii)), $\im(K_p(A)\to G^{\ab}(\dbZ_p))$ contains the Sylow $p$-subgroup of the image of the $m$-th power endomorphism of $G^{\ab}(\dbZ_p)$. So $\Ker( G^{\ab}(\dbZ_p)\to  G^{\ab}(\dbF_p))\leqslant\im(K_p(A)\to G^{\ab}(\dbZ_p))$. 

Conditions 2.4 (i), (iii), (iv) and (va) hold in the context of $G$ and $\im(K_p(A)\to G(\dbZ_p))$ (cf. 2.4.3 (a)). As $n+p>3$ (cf. (i)), the factors of the composition series of $G(k)$ are either cyclic of order prime to $p$ or non-abelian simple groups (cf. 2.2.4). So $K_p(A)$ surjects onto $\tilde G(\dbZ_p)$, cf. 2.4.2 (a). Due to (1), the group $\im(K_p(A)\to G(\dbF_p))$ surjects onto $G^{\ad}(\dbF_p)=\tilde G(\dbF_p)$. So $\im(K_p(A)\to G(\dbF_p))$ contains a normal subgroup $S$ generated by elements of order $p$ and surjecting onto the simple group $G^{\ad}(\dbF_p)^\prime$, cf. 2.2.4. So $S\leqslant G^{\der}(\dbF_p)$ and so $S=G^{\der}(\dbF_p)$, cf. 2.2.3. So $G^{\der}(\dbZ_p)\vartriangleleft \im(K_p(A)\to G(\dbZ_p))$, cf. 2.3. So the group $ G(\dbZ_p)/\im(K_p(A)\to G(\dbZ_p))$ is naturally a quotient of $G^{\ab}(\dbF_p)$. As $K_p(A)$ surjects onto $G^{\ad}(\dbF_p)=\tilde G(\dbF_p)$ and as the group $G^{\ad}(\dbF_p)/G^{\ad}(\dbF_p)^\prime$ has the same number of elements as $H^1(\dbF_p,G^{\ab}_{\dbF_p})$ and so as $G^{\ab}(\dbF_p)$, by reasons of prime to $p$ orders we get that the group $G(\dbZ_p)/\im(K_p(A)\to G(\dbZ_p))$ is trivial. So $K_p(A)$ surjects onto $G(\dbZ_p)$. 

As we have a short exact sequence $0\to T_1(\dbZ_p)\to H_{A\dbZ_p}(\dbZ_p)\to G(\dbZ_p)\to 0$, from (2) we get that $\scrV_p(A)^{-1}=[T_1(\dbZ_p):T_1(\dbZ_p)\cap K_p(A)]$. So we have $\scrV_p(A)^{-1}=up^s$, where $u\in\dbN$ divides $(p-1)$ and $s\in\dbN\cup\{0\}$.\endproof

\bigskip\noindent
{\bf 3.5. Variants of 3.4.} 
We assume $g.c.d.(p,c(\tilde G))>1$. So $\tilde G(\dbZ_p)\times \tilde G(\dbZ_p)$ has proper, closed subgroups surjecting onto each one of the two factors $\tilde G(\dbZ_p)$ as well as onto $\tilde G(\dbF_p)\times \tilde G(\dbF_p)$ (like the subgroup generated by $\im(\tilde G^1(\dbZ_p)\times \tilde G^1(\dbZ_p)\to \tilde G(\dbZ_p)\times \tilde G(\dbZ_p))$ and by the diagonal image of $\tilde G(\dbZ_p)$ in $\tilde G(\dbZ_p)\times \tilde G(\dbZ_p)$, where $\tilde G^1$ is an isogeny cover of $\tilde G$ of degree $p$). So the methods of 3.4 can be extended to cases in which $H_A^{\ad}$ is not a.s. only if by some reasons we can appeal (see 2.4 (iv) and 2.4.1) to ``sufficiently many" epimorphisms from $\Gal(E^{\text{conn}})$ onto groups of $\dbZ_p$-valued points of quotients of $H_{A\dbZ_p}^{\ab}$ or if we are in situations in which we can combine 2.2.5 and 2.4 with 3.4. Here is an example modeled on 3.4.1.1.

\medskip\noindent
{\bf 3.5.1. Corollary.} {\it We assume that $p=2$, that $T$ is a $\dbG_m$ group and that there is a semisimple subgroup $F$ of $H_{A\dbZ_2}^{\der}$ such that $c(F)$ is odd, the intersection $F\cap T$ is the identity section of $T$ and the simple factors of $F^{\ad}$ are not $PGL_2$, $PGL_3$, $PGSp_4$, $PGU_3$, $PGU_4$, or $\Res_{W(\dbF_4)/\dbZ_2} PGL_2$ groups. We also assume that $o(H_{A\dbZ_2}/(TF))$ is odd, that the adjoint of $H_{A\dbZ_2}^{\der}/F$ is a.s. of $D_{2n+3}$ or $A_{4n+3}$ Dynkin type, that $o(H_{A\dbZ_2}^{\der}/F)=2$, that $K_2(A)$ surjects onto $H_{A\dbZ_2}(\dbF_2)$ and that $E^{\text{conn}}$ is linearly disjoint from $\dbQ_{\mu_{2^{\infty}}}$. Then $K_2(A)$ surjects onto $(H_{A\dbZ_2}/T)(\dbZ_2)$.}  

\medskip
\proof
We have $F(\dbZ_2)\vartriangleleft K_2(A)$, cf. 2.3. So as in the proof of 3.4.1.1 (a), the Corollary follows from 2.4.1 (b) applied to $H_{A\dbZ_2}/F$ and $\im(K_2(A)\to (H_{A\dbZ_2}/F)(\dbZ_2))$.\endproof

\bigskip\noindent
{\bf 3.6. Remarks.} 
{\bf (a)} We refer to 3.4.1.1. Only its cases when $\tilde G$ either is of $A_{4n+1}$, $B_n$, $C_n$ or $D_{2n}$ Dynkin type, or is not adjoint of $D_{2n+3}$ Dynkin type, or is of $A_{4n+3}$ Dynkin type with $g.c.d.(2,o(\tilde G))=2$, are also implied by 3.3.1.

{\bf (b)} See [17] for Picard modular varieties. The simplest situations to which one can apply 3.4.1.1 (resp. 3.4.2.1) are related to Siegel (resp. Picard) modular varieties of dimension 6 (resp. 2) for $p=2$ (resp. $p=3$). Let now $p$ be arbitrary. In the case of Picard modular varieties, if $H_{A\dbZ_p}$ is a $GL_3\times\dbG_m$ group, then one can apply the classification of subgroups of $PGL_3(\dbF_p)$ (see [22]) to find general conditions under which (1) holds.

{\bf (c)} See [38] for many cases for which 3.1 is known and so to which 3.2 to 3.5 apply. Also we mention that except for the proof of 3.4.2.1, in 3.3 and 3.4.1 (resp. in all of 3.3 to 3.5) we can substitute the assumption that 3.1 holds for $(A,p)$ by a reference to Bogomolov theorem of [4] (resp. to the refinement of this theorem; see [38]) that asserts that always $Z(GL(V_p(A))$ (resp. $T_{\dbQ_p}$) is a subgroup of $G_p$. In connection to 3.3.2 we need to add that it can be checked that always $G_p$ contains all semisimple, normal subgroups of $H_{A\dbQ_p}^{\der}$ whose adjoints are of isotypic $A_1$ Dynkin type (for instance, see [38]). So almost always 1.1 (a) (i.e. 3.1) is implied by 1.1 (b) and (c).

\bigskip
\noindent
{\boldsectionfont \S 4. Examples}

\bigskip
Let $E$, $i_E$, $A$, $d$, $H_1$, $H_A$, $p$, $T_p(A)$, $H_{A\dbZ_p}$, $\rho_{A,p}$, $G_p$, $E^{\text{conn}}$ and $K_p(A)$ be as in the beginning of \S 1. In 4.1 and 4.3 we include two special cases of 3.4.1.1 (a) and (b) but stated in a simpler way (i.e. with much fewer hypotheses). In 4.1.1, 4.1.2 and 4.2 we include concrete examples pertaining to $p\Le 3$ and $d\in\{1,2,3\}$. In 4.4 we include remarks. 

Warning: until end we do not assume that conditions 1.1 (a) to (c) hold but only that $A$ has a polarization $p_A$ of degree prime to $p$. We denote also by $p_A$ the perfect, alternating form on $T_p(A)$ defined by $p_A$. For simplicity, we also denote 
$$GL_{2d}:=GL(T_p(A))\,\,\,\text{and}\,\,\,GSp_{2d}:=GSp(T_p(A),p_A).$$ 
As $p_A$ is defined over $E$, the group $\im(\rho_{A,p})$ is a subgroup of $GSp_{2d}(\dbZ_p)$. It is well known that if $E$ is linearly disjoint from $\dbQ_{\mu_{p^{\infty}}}$, then $\im(\rho_{A,p})$ surjects onto $GSp_{2d}^{\ab}(\dbZ_p)$ (the proof of this is the same as of the last part of 3.2.1). The following result is an extension to the case of all primes of a result of Serre for $p\Ge 5$ (see [33, pp. 50 to 53]). 

\bigskip\noindent
{\bf 4.1. Theorem.} {\it We assume that $E$ is linearly disjoint from $\dbQ_{\mu_{p^{\infty}}}$ and that $A$ has a polarization $p_A$ of degree prime to $p$. If $p=2$ we also assume that $d\Ge 3$ and if $p=3$ we also assume that $d\Ge 2$. Then $\im(\rho_{A,p})=GSp_{2d}(\dbZ_p)$ iff $\im(\rho_{A,p})$ surjects onto $PGSp_{2d}(\dbF_p)$.} 

\medskip
\proof 
The only if part is trivial. We check the if part. Let $K_1$ be the image of $\im(\rho_{A,p})$ in $GSp_{2d}(\dbF_p)$. The kernel of the epimorphism $K_1\to PGSp_{2d}(\dbF_p)$ has order prime to $p$ and the normal subgroup $PGSp_{2d}(\dbF_p)^\prime$ of $PGSp_{2d}(\dbF_p)$ is generated by elements of order $p$ (cf. 2.2.4) and has an index prime to $p$. So the subgroup $K_1^\prime$ of $K_1$ generated by elements of order $p$ surjects onto $PGSp_{2d}(\dbF_p)^\prime$. So as $K_1^\prime\leqslant Sp_{2d}(\dbF_p)$, we have $K_1=Sp_{2d}(\dbF_p)$ (cf. 2.2.3). So $Sp_{2d}(\dbZ_p)\vartriangleleft\im(\rho_{A,p})$, cf. 2.3. So $\im(\rho_{A,p})=GSp_{2d}(\dbZ_p)$.\endproof

\medskip\noindent
{\bf 4.1.1. A quartic curve with $p=2$.} Let $p=2$. Let $A$ be the Jacobian of the projective, smooth quartic curve $C$ over $\dbQ$ defined by the homogeneous equation
$$xz^3+zx^3+zx^2y+zy^3+x^4+x^3y+x^2y^2+y^4 =0.$$
So $d=3$ and $E=\dbQ$. Let $p_A$ be the canonical principal polarization of $A$. It is known that $\im(\rho_{A,2})$ surjects onto $GSp_6(\dbF_2)=Sp_6(\dbF_2)$, cf. [34, p. 69]. So from 4.1 we get that $$\im(\rho_{A,2})=GSp_6(\dbZ_2).$$
\noindent
{\bf 4.1.2. A hyperelliptic curve of genus $2$ with $p=3$.} Let $p=3$. Let $C$ be the smooth, projective curve over $\dbQ$ having $C_0:=\Spec(\dbQ[x,y]/(y^2-x^5+x-1))$ as an open, dense subscheme. It is a hyperelliptic curve of genus $d=2$. Let $A$ be the Jacobian of $C$ and let $p_A$ be its canonical principal polarization. It is known that $\im(\rho_{A,3})$ surjects onto $GSp_4(\dbF_3)$, cf. [11, p. 509]. So from 4.1 we get that $$\im(\rho_{A,3})=GSp_4(\dbZ_3).$$
\indent
Next we include examples that implicitly point out that the exclusions in 4.1 of the cases when $p=2$ and $d\in\{1,2\}$ are necessary.

\bigskip\noindent
{\bf 4.2. Hyperelliptic curves of genus at most $2$ with $p=2$.} 
Let $x$ and $y$ be independent variables. Let $d\in\{1,2\}$. Let $f_{3d}\in\dbQ[x]$ be a polynomial of degree $3d$ that has $3d$ distinct roots $\beta_1$, ..., $\beta_{3d}$ in $\overline{\dbQ}$. Let $F$ be the splitting field of $f_{3d}$. The Galois group $\Gal(F/\dbQ)$ is a subgroup of the symmetric group $S_{3d}$ of permutations of the set $\{\beta_1,...,\beta_{3d}\}$. Let $C$ be the smooth, projective curve over $\dbQ$ having $C_0:=\Spec(\dbQ[x,y]/(y^2-f(x)))$ as an open, dense subscheme. If $d=1$, then $C$ is an elliptic curve and if $d=2$, then  $C$ is a hyperelliptic curve of genus $d=2$. If $d=1$, let $P_1$ be the point of $C(\overline{\dbQ})\setminus C_0(\overline{\dbQ})$. If $d=2$, let $P_1$ and $P_2$ be the two distinct points of $C(\overline{\dbQ})\setminus C_0(\overline{\dbQ})$. All these points are defined over $\dbQ$. Let $A$ be the Jacobian of $C$. Let $\phi:C\hookrightarrow A$ be the embedding defined by the rule $P\to P-P_1$. Let $p_A$ be the canonical principal polarization of $A$. 

Let $p=2$. Let $X_i\in C_0(F)$ be of coordinates $(\beta_i,0)$. For $i$, $j\in\{1,...,3d\}$, with $i\neq j$, the point $X_{ij}:=\phi(X_i)-\phi(X_j)\in A(F)$ has order $2$. So $\Gal(F/\dbQ)$ normalizes the $\dbF_2$-subspace $\scrX$ of $T_2(A)/2T_2(A)=A_{\overline{\dbQ}}[2]$ generated by $X_{ij}$'s. 

\medskip\noindent
{\bf 4.2.1. Theorem.} {\it We have $\im(\rho_{A,2})=GSp_{2d}(\dbZ_2)$ iff the following two conditions hold

\medskip
{\bf (i)} the Galois group $\Gal(F/\dbQ)$ is the full symmetric group $S_{3d}$;

\smallskip
{\bf (ii)} the field $F$ does not contain any one of the following three fields $\dbQ(i)$, $\dbQ(\sqrt{2})$, $\dbQ(\sqrt{-2})$ (i.e. $F$ and $\dbQ(i,\sqrt{2})$ are linearly disjoint).}

\medskip
\proof
The group $GSp_{4}(\dbF_2)=Sp_{2d}(\dbF_2)$ is isomorphic to $S_{3d}$, cf. [2]. Let $K:=\im(\rho_{A,2})$ and let $K_0:=\Ker(Sp_{2d}(\dbZ_2)\to Sp_{2d}(\dbF_2))$. 

We first prove the ``only if" part. So we assume that $K=GSp_{2d}(\dbZ_2)$. The representation of $GSp_{2d}(\dbF_2)$ on $T_2(A)/2T_2(A)$ is irreducible and so we have $\scrX=T_2(A)/2T_2(A)$. So the image of $\im(\rho_{A,2})$ in $GSp_{2d}(\dbF_2)$ is naturally identified with a quotient of $\Gal(F/\dbQ)$ and so by reasons of orders we get that (i) holds. The group $GSp_{2d}(\dbZ_2)/K_0$ is isomorphic to $Sp_{2d}(\dbF_2)\times \dbG_m(\dbZ_2)=S_{3d}\times \dbG_m(\dbZ_2)$. Moreover we can choose these products such that $\Gal(\dbQ)$ acts on $\{\beta_1,...,\beta_{3d}\}$ via $S_{3d}=\Gal(F/\dbQ)$ and on $\dbZ_2$ via the $2$-adic cyclotomic character $\Gal(\dbQ)\to\dbG_m(\dbZ_2)$. So as (i) holds, $K$ surjects onto $GSp_{2d}(\dbZ_2)/K_0$ iff the field $F$ is linearly disjoint from $\dbQ_{\mu_{2^\infty}}$. The group $S_{3d}=\Gal(F/\dbQ)$ has only three normal subgroups: the trivial one, $A_{3d}$ and $S_{3d}$. So as (i) holds, we easily get that $F$ and $\dbQ_{\mu_{2^{\infty}}}$ are linearly disjoint over $\dbQ$ iff condition (ii) holds. 

We now check the ``if" part. So we assume that (i) and (ii) hold. We identify $\dbG_m=GSp_{2d}^{\ab}$. It is easy to see that the representation of $\Gal(F/\dbQ)$ on $\scrX$ is faithful. So as (i) holds, the image of $K$ in $GSp_{2d}(\dbF_2)$ has a subgroup isomorphic to $S_{2d}=\Gal(F/\dbQ)$. So by reasons of orders we get that $K$ surjects onto $GSp_{2d}(\dbF_2)$. Due to the iff's of the previous paragraph, as (ii) holds we get that $K$ surjects onto $GSp_{2d}(\dbZ_2)/K_0$. 

We first consider the case $d=2$. The group $K$ surjects onto $\dbG_m(\dbZ_2)$. So $K=GSp_{4}(\dbZ_2)$ iff $K^{\d}:=K\cap Sp_{4}(\dbZ_2)$ is $Sp_{4}(\dbZ_2)$. This is equivalent to $K^{\d}$ surjecting onto $Sp_{4}(\dbF_2)$ (cf. 2.2.5) and so to $K$ surjecting onto $GSp_{4}(\dbZ_2)/K_0$. So $K=GSp_{4}(\dbZ_2)$.

Let now $d=1$; so $GSp_2=GL_2$. Let $\grg:=\Lie(GL_{2\dbF_2})$. We view $\grh:=\Lie(PGL_{2\dbF_2})$ as a quotient of $\grg$. Let $\grk$ be the intersection of the image of $K$ in $GL_2(\dbZ_2/4\dbZ_2)$ with $\Ker(GL_2(\dbZ_2/4\dbZ_2)\to GL_2(\dbF_2))$. As $K$ surjects onto $GL_2(\dbF_2)$, $\grk$ is a $GL_2(\dbF_2)$-submodule of $\grg$. As $K$ surjects onto $\dbG_m(\dbZ_2)$, $\grk$ surjects onto $\grh/[\grh,\grh]$. It is well known that the short exact sequence of $GL_2(\dbF_2)$-modules $0\to [\grh,\grh]\to \grh\to \grh/[\grh,\grh]\to 0$ does not split (this is a particular case of [37, 3.10 1)]). So $\grk$ surjects onto $\grh$. Also it is well known that the short exact sequence of $GL_2(\dbF_2)$-modules $0\to \Lie(Z(GL_{2\dbF_2}))\to \grg\to \grh\to 0$ does not split. So we have $\grk=\grg=[\grg,\grg]$ and so $K$ surjects onto $GL_2(\dbZ_2/4\dbZ_2)$. From this and 2.2.1 (b) we get that $K$ surjects onto $GL_2(\dbZ_2/8\dbZ_2)$. So $K=GL_2(\dbZ_2)$, cf. 2.2.1 (a).\endproof

\bigskip\noindent
{\bf 4.3. Theorem.} {\it We assume that $p=2$, that $\End(A)\otimes_{\dbZ} \dbZ_2$ is an $M_2(\dbZ_2)$ matrix algebra, that $\End(A_{\overline{E}})\otimes_{\dbZ} \dbR$ is the standard $\dbR$-algebra of quaternions, that $E^{\text{conn}}=E$, that $E$ is linearly disjoint from $\dbQ_{\mu_{2^{\infty}}}$, and that $d=2n$, with $n\in 2\dbN+3$ such that $2n$ is not of the form $\binom {2^{m+1}}{2^m}$, with $m\in\dbN$. We also assume that $A$ has a polarization $p_A$ of odd degree and such that $\End(A)\otimes_{\dbZ} \dbZ_2$ as a $\dbZ_2$-subalgebra of $\End(T_2(A))$ is self dual with respect to $p_A$.  Let $C_{\dbQ_2}$ be the identity component of the subgroup of $GSp_{2d\dbQ_2}$ fixing all elements of $\End(A)\otimes_{\dbZ} \dbQ_2$. We also assume that the Zariski closure $C_{\dbZ_2}$ of $C_{\dbQ_2}$ in $GL_{2d}$ is a reductive group scheme over $\dbZ_2$. If the reduction mod $2$ of the factorization of $\rho_{A,2}$ through $C_{\dbZ_2}(\dbZ_2)$ is surjective, then $\im(\rho_{A,2})$ surjects onto $C_{\dbZ_2}^{\ad}(\dbZ_2)$.} 

\medskip
\proof
Let $B(\dbF)$ be the field of fractions of the Witt ring $W(\dbF)$ of $\dbF:=\overline{\dbF_2}$. It is well known that $C^{\ab}_{\dbZ_2}$ is a $\dbG_m$ group, that $C^{\ad}_{\dbZ_2}$ is a.s. of $D_n$ Dynkin type and that $o(C^{\der})=2$ (see [18, pp. 391 and 395] for the picture over $\dbC$). The representation of $C_{B(\dbF)}$ on $T_2(A)\otimes_{\dbZ_2} B(\dbF)$ is a direct sum of two irreducible $2n$ dimensional representations $V_1$ and $V_2$. A well known theorem of Faltings says that the representation of $G_{2B(\dbF)}$ on $V_i$ is irreducible (see [12, p. 9]; here $i\in\{1,2\}$). It is known that there is a finite family of cocharacters of the image $G_{2B(\dbF)}(i)$ of $G_{2B(\dbF)}$ in $GL(V_i)$ that act on $V_i$ via the trivial and the identical characters of $\dbG_m$ and such that $G_{2B(\dbF)}(i)$ is generated by $G_{2B(\dbF)}(i)(B(\dbF))$-conjugates of these cocharacters, cf. [25, 5.10]. The hypotheses on $n$ imply that for any $s\in\dbN\cup\{0\}$ we have $2n\neq\binom {2s}{s}$. So from [25, pp. 211 and 212] we get that $G_{2B(\dbF)}(i)$ is the image $C_{2B(\dbF)}(i)$ of $C_{2B(\dbF)}$ in $GL(V_i)$. So as $C_{2B(\dbF)}(i)$ is isomorphic to $C_{2B(\dbF)}$, we get $\dim_{\dbQ_2}(G_2)\Ge\dim_{\dbQ_2}(C_{2\dbQ_2})$. So as $G_2\leqslant H_{A\dbQ_2}\leqslant C_{\dbQ_2}$, by reasons of dimensions we get $G_2=H_{A\dbQ_2}=C_{\dbQ_2}$. So 3.1 holds and $H_{A\dbZ_2}=C_{\dbZ_2}$ is a reductive group scheme over $\dbZ_2$. So the Theorem follows from 3.4.1.1 (a).\endproof

\bigskip\noindent
{\bf 4.4. Remarks.} {\bf (a)} Due to Hilbert's irreducibility theorem (see [29, Th. 5 of p. 9 and Ch. 9 and 10]), for each surjectivity criterion of 3.3, 3.4, 3.5, 4.1, and 4.3 there are situations when all its hypotheses hold (cf. also 3.6 (c)). In particular, we get the existence of polynomials $f_{3d}$ as in 4.2 such that $F$ and $\dbQ(i,\sqrt{2})$ are linearly disjoint. Also we get that for any dimension $d\Ge 2$, there are many examples similar to the ones of 4.1.1 and 4.1.2 but over some number field $E$ (cf. also 4.1.1 and 4.1.2 for $d\in\{2,3\}$  and cf. the fact that the locus $\scrJ_d$ of Jacobians in the Siegel upper-half plane $\scrH_d$ is irreducible for $d\Ge 2$; for instance, see [1, \S4 of Ch. IV] for this last fact).

{\bf (b)} The proof of 4.2.1 also explains why in 2.3 we do need to exclude the $PGSp_4$ simple factors for $q=2$ (even if they are not excluded by 2.2.5). 

{\bf (c)} Theorem 4.2.1 implicitly classifies all elliptic curves over $\dbQ$ whose $2$-adic representation has image $GL_2(\dbZ_2)$. It is easy to see that [39, Th. 2.6] implies that if in 4.2 we have $d=2$ and 4.2.1 (i) holds, then $H_{A\dbZ_2}=GSp_{4}$; so $E^{\text{conn}}=E=\dbQ$ and so the group $K$ of the proof of 4.2.1 is $K_2(A)$ itself. Also in 4.1.1 and 4.1.2 we have $H_{A\dbZ_p}=GSp_{2d}$, $E^{\text{conn}}=E$ and $K_p(A)=\im(\rho_{A,p})=GSp_{2d}(\dbZ_p)$ and so $\scrV_p(A)=1$.

{\bf (d)} We refer to 4.3. The subgroup $C_{1\dbZ_2}$ of $GL_{2d}$ fixing $\End(A)\otimes_{\dbZ} \dbZ_2$ is a $GL_{2n}$ group. The group scheme $C_{\dbZ_2}$ is reductive iff the standard involution of $C_{1\dbZ_2}$ defined by $p_A$ is mod $2$ alternating (see [5, 23.5 and 23.6] for the picture over $\overline{\dbF_2}$).

\bigskip
\noindent
{\it Acknowledgements.} We thank U of Arizona for good conditions for the writing of this Part II, the referee for many useful comments, and Serre and Dieulefait for mentioning to us the curves of 4.1.1 and respectively of 4.1.2.

\bigskip
\noindent
\references{37}
{\nspace{

\bigskip
\Ref[1]
E. Arbarello, M. Cornalba, P. A. Griffiths, and J. Harris,
\sl Geometry of algebraic curves,
\rm Grund. math. Wiss., Vol. {\bf 267}, Springer-Verlag, 1985.

\Ref[2]
J. H. Conway, R. T. Curtis, S. P. Norton, R. A. Parker, and R. A. Wilson,
\sl Atlas of finite groups,
\rm xxxiv+252 p, Oxford Univ. Press, Eynsham, 1985.

\Ref[3]
S. Bosch, W. L\"utkebohmert, and M. Raynaud,
\sl N\'eron models,
\rm Springer-Verlag, 1990.

\Ref[4]
F. A. Bogomolov,
\sl Sur l'alg\'ebricit\'e des repr\'esentations $l$-adiques,
\rm C. R. Acad. Sc. Paris Ser I Math. {\bf 290} (1980), pp. 701--703.

\Ref[5]
A. Borel,
\sl Linear algebraic groups, 
\rm Grad. Texts in Math., Vol. {\bf 126}, Springer-Verlag, 1991.

\Ref[6]
N. Bourbaki,
\sl Lie groups and Lie algebras, 
\rm Chapters {\bf 4--6}, Springer-Verlag, 2002.

\Ref[7]
F. Bruhat and J. Tits, 
\sl Groupes r\'eductifs sur un corps local. II.,
\rm Inst. Hautes \'Etudes Sci. Publ. Math., Vol. {\bf 60}, pp. 5--184, 1984.

\Ref[8]
P. Deligne,
\sl Vari\'et\'es de Shimura: Interpr\'etation modulaire, et
techniques de construction de mod\`eles canoniques,
\rm Proc. Sympos. Pure Math., Vol. {\bf 33}, Part 2, pp. 247--290, 1979. 

\Ref[9]
P. Deligne,
\sl Hodge cycles on abelian varieties,
\rm Hodge cycles, motives, and Shimura varieties, Lecture Notes in Math., Vol. {\bf 900}, pp. 9--100, Springer-Verlag, 1982.

\Ref[10]
M. Demazure, A. Grothendieck, et al. 
\sl Sch\'emas en groupes, Vol. III,
\rm Lecture Notes in Math., Vol. {\bf 153}, Springer-Verlag, 1970.

\Ref[11]
L. V. Dieulefait,
\sl Explicit determination of the images of the Galois representations attached to abelian surfaces with $\End(A)=\dbZ$,
\rm Experimental Math. {\bf 11}, No. 4 (2002), pp. 503--512.

\Ref[12]
G. Faltings,
\sl Finiteness theorems for abelian varieties over number fields,
\rm Arithmetic Geometry, pp. 9--27 (english translation of Inv. Math. {\bf 73}, No. 3 (1983), pp. 349--366), Springer-Verlag, 1986.

\Ref[13] 
J.-M. Fontaine, 
\sl Le corps des p\'eriodes p-adiques, 
\rm J. Ast\'erisque {\bf 223}, pp. 59--101, Soc. Math. de France, Paris, 1994.

\Ref[14] 
J.-M. Fontaine, 
\sl Repr\'esentations p-adiques semi-stables,
\rm J. Ast\'erisque {\bf 223}, pp. 113--184, Soc. Math. de France, Paris, 1994.

\Ref[15]
D. Gorenstein, R. Lyons, and R. Soloman,
\sl The classification of the finite simple groups, Number 3,
\rm Math. Surv. and Monog., Vol. {\bf 40}, No. 3, Amer. Math. Soc., 1994.

\Ref[16]
J. E. Humphreys, 
\sl Conjugacy classes in semisimple algebraic groups, 
\rm Math. Surv. and Monog., Vol. {\bf 43}, Amer. Math. Soc., 1995.

\Ref[17] B. B. Gordon, 
\sl Canonical models of Picard modular surfaces, 
\rm The Zeta function of Picard modular surfaces, pp. 1--28, Univ. Montr\'eal Press, Montreal, QC, 1992.

\Ref[18]
R. E. Kottwitz, 
\sl Points on some Shimura varieties over finite fields, 
\rm J. of Amer. Math. Soc. {\bf 5}, No. 2 (1992), pp. 373--444.

\Ref[19]
M. Larsen,
\sl Maximality of Galois actions for compatible systems,
\rm Duke Math. J. {\bf 80}, No. 3 (1995), pp. 601--630. 

\Ref[20]
M. Larsen and R. Pink,
\sl A connectedness criterion for $l$-adic Galois representations,
\rm Israel J. Math. {\bf 97} (1997), pp. 1--10.

\Ref[21]
B. Mazur,
\sl Modular curves and the Eisenstein ideal,
\rm Inst. Hautes \'Etudes Sci. Publ. Math., Vol. {\bf 47}, pp. 33--186, 1977.

\Ref[22]
H. Mitchell, 
\sl Determination of the ordinary and modular ternary linear groups,
\rm Trans. Amer. Math. Soc. {\bf 12} (1911), pp. 207--242.

\Ref[23]
D. Mumford, 
\sl Families of abelian varieties,
\rm Algebraic Groups and Discontinuous Subgroups, Proc. Sympos. Pure Math., Vol. {\bf 9}, pp. 347--352, Amer. Math. Soc., Providence, R. I., 1966.

\Ref[24]
R. Noot, 
\sl Lifting Galois representations and a conjecture of Fontaine and Mazur,
\rm Documenta Math. {\bf 6} (2001), pp. 419--445.

\Ref[25]
R. Pink,
\sl $l$-adic algebraic monodromy groups, cocharacters, and the Mumford--Tate conjecture,
\rm J. reine angew. Math. {\bf 495} (1998), pp. 187--237.

\Ref[26]
I. Satake,
\sl Holomorphic imbeddings of symmetric domains into a Siegel space,
\rm  Amer. J. Math. {\bf 87} (1965), pp. 425--461.

\Ref[27] 
J.-P. Serre,
\sl Propri\'et\'es galoisiennes des points d'ordre fini des courbes elliptiques,
\rm Invent. Math. {\bf 15}, No. 4 (1972), p 259--331. 

\Ref[28]
J.-P. Serre,
\sl Abelian $l$-adic representations and elliptic curves,
\rm Addison--Wesley Publ. Co., Redding, Mass., 1989. 

\Ref[29] 
J.-P. Serre, 
\sl Lectures on the Mordell--Weil theorem,
\rm Asp. of Math., Vol. {\bf E15},  Friedr. Vieweg and Sohn, Braunschweig, 1989.

\Ref[30]
J.-P. Serre, 
\sl Propri\'et\'es conjecturales des groupes de Galois motiviques et des repr\'esentations $l$-adiques,
\rm Proc. Sympos. Pure Math., Vol. {\bf 55}, Part 1, pp. 377--400, Amer. Math. Soc., 1994.

\Ref[31]
J.-P. Serre,
\sl Travaux de Wiles (et Taylor, ...), Partie I,
\rm S\'em. Bourbaki, Vol. 1994/95, Exp. No. 803, 5, pp. 319--332, Ast\'erisque {\bf 237}, Soc. Math. de France, Paris, 1996.

\Ref[32]
 J.-P. Serre, 
\sl Galois Cohomology, 
\rm Springer-Verlag, 1997.

\Ref[33]
J.-P. Serre, 
\rm Collected papers, Vol. IV, Springer-Verlag, 2000.

\Ref[34]
T. Shioda,
\sl Plane quartics and Mordell--Weil lattices of type  E7,
\rm Comm. Math. Univ. Sancti Pauli {\bf 42}, No. 1 (1993), pp. 61--79.

\Ref[35]
J. Tits,
\sl Reductive groups over local fields, 
\rm Proc. Sympos. Pure Math., Vol. {\bf 33}, Part 1, pp. 29--69, Amer. Math. Soc., 1979.

\Ref[36] 
T. Tsuji, 
\sl p-adic \'etale cohomology and crystalline cohomology
in the semi-stable reduction case, 
\rm Invent. Math. {\bf 137} (1999), pp. 233--411.

\Ref[37]
A. Vasiu,
\sl Surjectivity criteria for p-adic representations, Part I,
\rm Manuscripta Math. {\bf 112}, No. 3 (2003), pp. 325--355.

\Ref[38]
A. Vasiu,
\sl The Mumford--Tate conjecture and Shimura varieties, Part I,
\rm manuscript Sep. 2003 (available at xxx.lanl.gov).

\Ref[39]
Y. Zarhin,
\sl Very simple $2$-adic representations and hyperelliptic Jacobians,
\rm Moscow Math. J. {\bf 2}, No. 2 (2000), pp. 403--431.

\Ref[40]
W. Waterhouse, 
\sl Introduction to affine group schemes, 
\rm Grad. Texts in Math., Vol. {\bf 89}, Springer-Verlag, 1974.

\Ref[41] J.-P. Wintenberger, 
\sl  Une extension de la th\'eorie de la multiplication complexe,
\rm J. reine angew. Math. {\bf 552} (2002), pp. 1--14. 

\Ref[42] J.-P. Wintenberger, 
\sl D\'emonstration d'une conjecture de Lang dans des cas particuliers,
\rm  J. reine angew. Math. {\bf 553} (2002), pp. 1--16.

}}
\enddocument